\documentclass[notitlepage,leqno,10pt]{article}

\usepackage{latexsym}
\usepackage{amsmath}
\usepackage{amsfonts}
\usepackage{amssymb}
\usepackage{hyperref}
\renewcommand{\a }{\alpha }
\renewcommand{\b }{\beta }
\renewcommand{\d}{\delta }
\newcommand{\D }{\Delta }

\newcommand{\e }{\varepsilon }
\newcommand{\g }{\gamma}

\newcommand{\G }{\Gamma }
\renewcommand{\l }{\lambda }

\newcommand{\n }{\nabla }
\newcommand{\var }{\varphi }

\renewcommand{\k }{\kappa }

\renewcommand{\O }{\Omega }
\newcommand{\z }{\zeta }
\newcommand{\ov}{\overline}
\newcommand{\no}{\noindent}
\newcommand{\ms}{\medskip}
\newcommand{\pa}{\partial}
\newcommand{\tr }{\hbox{ tr } }
\newcommand{\intbar}{\mathop{\int\makebox(-13.5,0){\rule[4pt]{.7em}{0.3pt}}%
\kern-6pt}\nolimits}

\newcommand{\be}{\begin{equation}}
\newcommand{\ee}{\end{equation}}
\newenvironment{pf}{\noindent{\sc Proof}.\enspace}{\rule{2mm}{2mm}\medskip}
\newenvironment{pfn}{\noindent{\sc Proof}}{\rule{2mm}{2mm}\medskip}

\newtheorem{corollary}{Corollary}[section]

\newcommand{\R}{\mathbb{R}}

\newcommand{\N}{\mathbb{N}}

\author{Fethi Mahmoudi $^1$,  Andrea Malchiodi $^2$ and Juncheng Wei $^3$}

\date{}

\title{Transition Layer for the Heterogeneous Allen-Cahn Equation}

\begin{document}

\newtheorem{lem}{Lemma}[section]
\newtheorem{pro}[lem]{Proposition}
\newtheorem{thm}[lem]{Theorem}
\newtheorem{rem}[lem]{Remark}
\newtheorem{cor}[lem]{Corollary}
\newtheorem{df}[lem]{Definition}

\maketitle

\begin{center}

$^{1,2} $ {\small SISSA, via Beirut 2-4, 34014 Trieste, Italy.}
\\

$^3$ {\small Department of Mathematics,
  The Chinese University of Hong Kong,
  Shatin, Hong Kong.}

\end{center}

\footnotetext[1]{E-mail addresses: mahmoudi@ssissa.it
(F.Mahmoudi),} \footnotetext[2]{E-mail addresses:
malchiod@sissa.it (A. Malchiodi),} \footnotetext[3]{E-mail
addresses: wei@math.cuhk.edu.hk (J. Wei)}

\

\

\noindent {\sc abstract}.  We consider the equation
\begin{equation}\label{eq:fife}
\e^{2}\Delta u =(u-a(x))(u^2-1)\ \  \mbox{in} \ \Omega, \ \
 \frac{ \partial  u}{\partial \nu} =0 \ \ \mbox{on}\ \partial \Omega,
\end{equation}
where $\Omega $ is a smooth and bounded domain in $\R^n$, $\nu$ the
outer unit normal to $\pa\O$, and $a$ a smooth function satisfying
$-1<a(x)<1$ in $\ov{\O}$. We set $K$, $\O_+$ and $\O_-$ to be
respectively the zero-level set of $a$, $\{a>0\}$ and $\{a<0\}$.
Assuming $\n a \neq 0$ on $K$ and $a\ne 0$ on $\pa \O$, we show that
there exists a sequence $\e_j \to 0$ such that equation
\eqref{eq:fife} has a solution $u_{\e_j}$ which converges uniformly
to $\pm 1$ on the compact sets of $\O_{\pm}$ as $j \to + \infty$.
This result settles in general dimension a conjecture posed in
\cite{fg}, proved in \cite{dkw} only for $n=2$.

\begin{center}

\bigskip\bigskip

\noindent{\it Key Words:} Fife-Greenlee problem, heterogeneous
Allen-Cahn equation, Interior transition layers, spectral gaps.

\bigskip

\centerline{\bf AMS subject classification: 35J25, 35J40, 35B34,
35B40}

\end{center}

\section{Introduction}\label{s:in}

Given a smooth bounded domain $\O$ of $\R^n$ ($n\ge 2$), we consider
the following problem
\begin{equation}\label{eq:fife1}
    \left\{
      \begin{array}{ll}
        \e^2 \D u = h(x,u) & \hbox{ in } \O, \\
        \frac{\pa u}{\pa \nu} = 0 & \hbox{ on } \pa \O,
      \end{array}
    \right.
\end{equation}
where $\e$ is a small parameter, $\nu$ the unit outer normal vector
to $\pa \O$ and $h$ a smooth function such that  the equation
$h(x,t)=0$ admits two different stable solutions $t_1\ne t_2$ for
any $x\in\ov\O$. Using matched asymptotics,  Fife and Greenlee in
\cite{fg} proved under some hypothesis on $h$ the existence of a
solution of \eqref{eq:fife1} which converges uniformly to $t_i$ in
the compact subsets of $\O_i,\; i=1,2$, where $\O_1$ and $\O_2$ are
two subdomains of $\O$ such that $\ov\O={\ov\O_1}\cup {\ov\O_2}$.

\

\no In this paper we consider a {\em model} heterogeneous case
$h(x,u) = (u-a(x))(u^2-1)$, for a smooth function $a$ satisfying
$-1<a(x)<1$ on $\ov{\O}$ and $\n a \neq 0$ on the set $K=\{
a(x)=0\}$, with $K \cap \pa \O = \emptyset$. We prove the existence
of a new type of solution of \eqref{eq:fife1} for any $n\ge 2$
settling in full generality a result previously proved in \cite{dkw}
for the particular case $n=2$.

\

Let us describe the result in more detail: in the case
$h(x,u)=(u-a(x))(u^2-1)$ problem \eqref{eq:fife1} becomes
\begin{equation}\label{eq:fife2}
    \left\{
      \begin{array}{ll}
        \e^2 \D u = (u-a(x))(u^2-1) & \hbox{ in } \O, \\
        \frac{\pa u}{\pa \nu} = 0 & \hbox{ on } \pa \O.
      \end{array}
    \right.
\end{equation}
In particular, when $a\equiv 0$, \eqref{eq:fife2} is nothing but the
Allen-Cahn equation  in material sciences (see \cite{ac})
\begin{equation}
\label{eq:ak} \left\{
\begin{array}{ll}
\e^{2}\Delta u + u-u^3 = 0 &  \mbox{in} \ \Omega, \\
 \frac{ \partial  u}{\partial \nu} =0 & \mbox{on}\ \partial \Omega.
\end{array}
\right.
\end{equation}
Here the function $u(x)$ represents a continuous realization of the
phase present in a material confined to the region $\Omega$ at the
point $x$.  Of particular interest are the solutions  which, except
for a narrow region, take values close to $+1$ or $-1$. Such
solutions are called {\em transition layers}, and have been studied
by many authors, see for instance \cite{acf, bs, fp, ks, k, mnw, mw,
mo, na, nt, prit, pt, r1, r2, sz}, and the references therein for
these and related issues.

In this paper, we are interested in transition layers for  the {\it
heterogeneous} equation (\ref{eq:fife2}). Define
$$
  K = \left\{ x \in \O \; : \; a(x) = 0 \right\}.
$$
We assume that $K$ is a smooth closed hypersurface  of $\O$ which
separates the domain into two disjoint components
\begin{equation}
\O=\O_{-} \cup K \cup \O_{+},
\end{equation}
with
\begin{equation}
 a(x) < 0 \,\hbox{  in  }  \Omega_{-} ,\qquad  a(x) > 0\,
\hbox{ in }\Omega_{+},\qquad  \n a \neq 0 \, \hbox{ on } K.
\label{conds}\end{equation}

 We then define the  Euler
functional $J_\e(u)$ associated
 to \eqref{eq:fife2} in $\O$  as
\begin{equation}\label{eq:Euler0}
 J_\e(u)=\frac{\e^2}{2}\int_\O |\n u|^2+\int_\O F(x,u)dx,
\end{equation}
where
\begin{equation*}
F(x,u):=\int_{-1}^u (s-a(x))(s^2-1)\,ds.
\end{equation*}
The solution constructed by Fife and Greenlee in \cite{fg} (adapted
to our choice of the function $h$) consists in adding an interior
transition layer correction to expressions of the form $t_i+\e
t_i^1+\e^2 t_i^2$, which approximate the solution $u$ in the regions
$\O_i$ (notice that with our choice of the function $h$, we have
$\O_1=\O_+$, $\O_2=\O_-$, $t_1\equiv -1$ and $t_2\equiv 1$). This
allowed Fife and Greenlee to construct an approximation $U_\e$ which
yields an exact solution of \eqref{eq:sfife2} using a classical
implicit function argument. No restriction on $\e$ are required, and
the solution satisfies
\begin{equation}\label{eq:asas}
u_\e\to -1\quad \hbox{in } \O_+ \quad \hbox{and }\quad u_\e\to
1\quad \hbox{in } \O_-
 \quad \hbox{as }\e\to 0.
\end{equation}

Super-subsolutions were later used by Angenent, Mallet-Paret and
Peletier in the one dimensional case (see \cite{amp}) for
construction and classification of stable solutions. Radial
solutions  were found variationally by Alikakos and Simpson in
\cite{as}. These results were extended by del Pino in \cite{d1} for
general (even non smooth) interfaces in any dimension, and further
constructions have been done recently by Dancer and Yan \cite{dy1}
and  Do Nascimento \cite{donascimento}. In particular, it was proved
in \cite{dy1} that solutions with the asymptotic behavior like
\eqref{eq:asas} are typically  minimizer of $J_\e$. Related results
can be found in \cite{ab1,abc1}.

\medskip
On the other hand, a solution exhibiting a transition layer in the
{\em opposite direction}, namely
\begin{equation}\hbox{ $u_\e \to +1$ in
$\Omega_{+}$\, , \quad  \ $u_\e \to -1$ on $\Omega_{-}$\  } \quad
\hbox{as }\e\to 0 \label{layer1}
\end{equation}
has been  believed to exist for many years. Hale and Sakamoto
\cite{hs}  established the existence of this type of solution in the
one-dimensional case, while this was done for the radial case in
\cite{d2}, see also \cite{dy}. The layer with the asymptotics in
\eqref{layer1} in this scalar problem is meaningful in describing
pattern-formation for reaction-diffusion systems such as
Gierer-Meinhardt with saturation, see
\cite{d2,f1,nishiura,sakamoto,s2} and the references therein.

For one-dimensional or radial problems it is possible to use
finite-dimensional reductions, which basically consist in
determining the location of the transition layer. In this kind of
approach, the same technique works for both the asymptotic behaviors
in \eqref{eq:asas} and \eqref{layer1}: the only difference is the
sign of the small eigenvalue (of order $\e$) arising from the
approximate degeneracy of the equation (when we tilt the solutions
perpendicularly to the interface). This makes the former solution
stable and the latter unstable.

On the other hand, one faces  a dramatically different situation in
 higher-dimensional, non-symmetric cases.  This is clearly
seen already linearizing  around a spherically symmetric solution of
\eqref{eq:fife} (with profile as in \eqref{layer1}), as bifurcations
of non-radial solutions along certain infinite discrete set of
values for $\e\to 0 $ take place, as established by Sakamoto in
\cite{s2}. This reveals that the radial solution has Morse index
which changes with $\e$ (precisely diverges as $\e \to 0$, as shown
in  \cite{dn}). This poses a serious difficulty for a general
construction. A phenomenon of this type was previously observed in
the one-dimensional case by Alikakos, Bates and Fusco \cite{abf1} in
finding solutions with any prescribed Morse index.

In \cite{dkw}, del Pino, Kowalczyk and the third author considered
the two-dimensional case, constructing transition layer solutions
with asymptotics as in \eqref{layer1}, while in this paper we extend
that result to any dimension.   Our main theorem is the following.

\vskip 0.5cm

\begin{thm}\label{t:mmm} Let $\O$ be a smooth bounded domain of $\R^n$
($n\ge 2$) and assume that $a : \ov{\O} \to (-1,1)$ is a smooth
function. Define $K$, $\O_+$ and $\O_-$ to be respectively the
zero-set, the positive set and the negative set of $a$. Assume that
$\n a \neq 0$ on $K$ and that $K \cap \pa \O = \emptyset$. Then
there exists a sequence $\e_j \to 0$ such that problem
\eqref{eq:fife2} has a solution $ u_{\e_j}$ which approach $1$ in
$\O_+$ and $-1$ in $\O_-$. Precisely, parameterizing a point $x$
near $K$ by $x=(\ov{y},\ov{\z})$, with $\ov{y} \in K$ and
$\ov{\z}=d(x, K)$ (with sign, positive in $\O_+$), $u_{\e_j}$ admits
the following behavior
\begin{equation*}
  u_{\e_j}(\ov{y},\ov{\z})= H\left( \frac{\ov{\z}}{\e_j} +\Phi(\ov{y})
  \right)+O(\e_j) \quad \hbox{as } j\to +\infty.
\end{equation*}
Here $\Phi$ is a smooth function defined on $K$ and $H(\z)$ is the
unique hetheroclinic solution of
\begin{equation}
\label{ode} H^{''} +H -H^3=0, \qquad \quad  H(0)=0, \qquad \quad
H(\pm \infty) = \pm 1.
\end{equation}
\end{thm}

\no As  in \cite{dkw0}, \cite{dkw}, \cite{mm4}-\cite{mm1}, \cite{mw}
and other results for singularly perturbed (or geometric) problems,
the existence is proved only along a sequence $\e_j \to 0$ (actually
 it can be obtained for $\e$ in a sequence of
intervals $(a_j, b_j)$ approaching zero, but not for any small
$\e$). This is caused by a resonance phenomenon we are going to
discuss below, explaining the ideas of the proof.

To describe the reasons which causes the main difficulty in proving
Theorem \ref{t:mmm}, we first scale  problem \eqref{eq:fife2} using
the change of variable $x\mapsto \e x$, so equation \eqref{eq:fife2}
becomes
\begin{equation}\label{eq:sfife2}
    \left\{
      \begin{array}{ll}
         \D u = (u-a(\e x))(u^2-1) & \hbox{ in } \O_\e, \\
        \frac{\pa u}{\pa \nu} = 0 & \hbox{ on } \pa \O_\e,
      \end{array}
    \right.
\end{equation}
where $\O_\e = \frac 1 \e \O$. Near the hypersurface $ K_\e := \frac
1 \e K$, we can choose scaled coordinates $(y,\z)$ in $\O_\e$ with
$y\in K_\e$ and $\z=dist(x,K_\e)$ (with sign), see Subsection
\ref{ss:gb}, and we let $\tilde{u}_\e$ denote the scaling of $u_\e$
to $\O_\e$: with these notations  we have that $\tilde{u}_\e (y,\z)
= u_\e(y, \e \z) \simeq H(\z)$. The function
$H(\z)=H\left(dist(x,K_\e)\right)$ for $x \in \O_\e$ can be then
considered as a first order {\em approximate solution} to
\eqref{eq:sfife2}, so it is natural to use local inversion arguments
near this function in order to find true solutions. For this purpose
it is necessary to understand the spectrum of the linearization of
\eqref{eq:sfife2} at approximate solutions.

\ Letting $L_\e$ be the linearization of \eqref{eq:sfife2} at
$\tilde{u}_\e$, it turns out that $L_\e$ admits a sequence of small
positive eigenvalues of order $\e$. Using asymptotic expansions (see
Section \ref{s:aa}, and in particular formula \eqref{eq:lje}), one
can see that this family behaves qualitatively like $\e - \e^2
\l_j$, where the $\l_j$'s are the eigenvalues of the
Laplace-Beltrami operator of $K$. By the Weyl's asymptotic formula,
we have that $ \l_j \simeq j^{\frac{2}{n-1}}$ as $j \to + \infty$,
therefore we have an increasing number of positive eigenvalues, many
of which accumulate to zero and sometimes, depending on the value of
$\e$, we even have the presence of a kernel: this clearly causes
difficulties if one wants to apply local inversion arguments. Notice
that, by the above qualitative formula, the average spectral gap of
resonant eigenvalues is of order $\e^{\frac{n+1}{2}}$. For the case
$n = 2$ (considered in \cite{dkw}) this gap is relatively large, so
it was possible to show invertibility using direct estimates on the
eigenvalues. However in higher dimension this is not possible
anymore, and one needs to apply different arguments.

\

To overcome this problem, we use an approach introduced in
\cite{mm1}, \cite{mm2} (see also \cite{mahm}, \cite{mm4},
\cite{mm3}) to handle similar resonance phenomena for another class
of singularly perturbed equations. The main idea consists in looking
at the eigenvalues (of the linearized problem) as functions of the
parameter $\e$, and  estimate their derivatives with respect to
$\e$. This can be rigorously done employing a classical theorem due
to T.Kato, see Proposition \ref{p:epS}, and by characterizing the
eigenfunctions corresponding
 to resonant modes. Using this result we get invertibility along a
suitable sequence $\e_j \to 0$, and the norm of the inverse operator
along this sequence has an upper bound of order
$\e_j^{-\frac{n+1}{2}}$ (consistently with the above heuristic
evaluation of the spectral gaps). This loss of uniform bounds as $j
\to + \infty$ should be expected, since more and more eigenvalues
are accumulating near zero. However, we are able to deal with this
further difficulty by choosing approximate solutions with a
sufficiently high accuracy.

\

Fixing an integer $k\ge 1$ and using the coordinates introduced
after \eqref{eq:sfife2}, from the fact that $a$ vanishes on $K$, one
can consider the Taylor expansion
$$a(\e y,\e \z)=\e \z b(\e y)+\sum_{l=2}^k (\e \z)^l b_l(\e y)+\tilde{b}(y,\z)
\quad\hbox{with}\quad \left| \tilde{b}(y,\z) \right|\le C_k
|\e\z|^{k+1},
$$
and look at an approximate solution of the form
$$
 u_{k,\e}(y,\z)=H(\z-\Phi(\e y))+\sum_{i=1}^k \e^i h_i(\e y,\z-\Phi(\e y)),
$$
for a smooth function $\Phi(\e y)=\Phi_0(\e
y)+\sum_{i=1}^{k-1}\e^i\Phi_i(\e y)$ defined on $K$ and some
corrections $h_i$ defined on $K\times \R_+$. Using similar Taylor
expansions of the Laplace-Beltrami operator in the above
coordinates, see Subsection \ref{ss:gb}, the couple
$(h_{j},\Phi_{j-1})$ for $j\ge 1$ can be determined via equations of
the form
\begin{equation}\label{eq:hN0}
\left\{
  \begin{array}{ll}
    \mathcal{L}_0 h_1= & -\kappa(\e y) H'(s)+(s+\Phi_0) b(\e y) (1-H^2(s)) \\
    \mathcal{L}_0 h_j=& \Phi_{j-1} b(\e y) (1-H^2(s))+\mathfrak{F}_k(s,\Phi_0,
    \dots,\Phi_{j-2},h_1,\dots,h_{j-1},b_1,\dots,b_j),\hbox{ for } j\ge
    2,
  \end{array}
\right.
\end{equation}
where $ \mathcal{L}_0 u = u'' + (1-3H^2) u$, $\mathfrak{F}_k$ is a
smooth function on its argument, and $s = \z - \Phi(\e y)$.
\eqref{eq:hN0} is always solvable in $h_j$ by the Fredholm
alternative if we choose properly the functions $\Phi_l$.

Such an accurate approximate solution allow us, using the above
characterization of the spectrum of the linearized operator and the
bound on its inverse, to apply the contraction mapping theorem to
find true solutions. Specifically for the homogeneous Allen-Cahn
equation, a related method was used in \cite{mw} to study the effect
of $\pa \O$ on the structure of solutions to \eqref{eq:ak}. Some
common arguments are here simplified, and we believe our approach
could also be used to handle general nonlinearities as in \cite{fg}.

\

\no The paper is organized in the following way: in Section
\ref{s:not} we collect some preliminary results concerning the
profile $H$, we expand the Euclidean metric and the
 Laplace-Beltrami operator in suitable coordinates near $K_\e$,
and recall some well-known spectral results. In Section \ref{s:aa}
we first construct approximate solutions, and then derive some
spectral properties of the linearized operator characterizing the
resonant eigenfunctions: this is a crucial step to apply Kato's
theorem.  Finally, Section \ref{s:pf} is devoted to the proof of our
main result.

\section{Notation and preliminaries}\label{s:not}

In this section we first collect some notation and conventions.
Then, we list some properties of the hetheroclinic solution $H$, and
we expand the metric and the Laplace-Beltrami operator in a local
normal coordinates. Finally we recall some results in spectral
theory like the Weyl asymptotic formula.

\begin{center}
\bf  Notation and convention
\end{center}
\no We shall always use the convention that capital letters like $A,
B, \dots$ will vary between $1$ and $n$, while indices like $i, j,
\dots$ will run between $1$ and $n-1$. We adopt the standard
geometric convention of summing over repeated indices.

\ms \no $y_{1}, \dots, y_{n-1}$ will denote coordinates in
$\R^{n-1}$, and they will also be written as $y=(y_{1}, \dots,
y_{n-1})$, while coordinates in $\R^n$ will be written
$x=(y,\z)\in \R^{n-1}\times\R$.

\ms \no  The hypersurface $K$ will be parameterized with local
coordinates $\ov{y} = (\ov{y}_1, \dots, \ov{y}_{n-1})$. It will be
convenient to define its dilation $K_\e := \frac 1 \e K$ which will
be parameterized by coordinates $(y_1, \dots, y_{n-1})$ related to
the $\ov{y}$'s simply by $\ov{y} = \e y$.

\ms \no  Derivatives with respect to the variables $\ov{y}$, $y$ or
$\z$ will be denoted by $\pa_{\ov{y}}$, $\pa_y$, $\pa_\z$ and for
brevity we shall sometimes use the notations $\pa_i$ for
$\pa_{y_i}$. When dealing with functions depending just on the
variable $\z$ we will write $H', \,h',\,\cdots$ instead of $\pa_\z
H,\,\pa_\z h,\,\cdots$.

\ms \no In a local system of coordinates, $(\ov{g}_{ij})_{i j}$
are the components of the metric on $K$ naturally induced by
$\R^{n-1}$. Similarly, $(\ov{g}_{AB})_{AB}$ are the entries of the
metric on $\O$ in a neighborhood of the hypersurface $K$.
$(\kappa_{i}^j)_{ij}$ will denote the components of the mean
curvature operator of $K$ into $\R^{n-1}$.

\

\no For a real positive variable $r$ and an integer $m$, $O(r^m)$
(resp. $o(r^m)$) will denote a function for which $\left|
\frac{O(r^m)}{r^m} \right|$ remains bounded (resp. $\left|
\frac{o(r^m)}{r^m} \right|$ tends to zero) when $r$ tends to zero.
For brevity, we might also write $O(1)$ (resp. $o(1)$) for a
quantity which stays bounded (resp. tends to zero) as $\e$ tends to
zero.

\subsection{Some analytic properties of the hetheroclinic solution  $H$}\label{ss:H}
In this subsection we collect some useful properties of the
hetheroclinic solution $H$ to \eqref{ode}. Note first that $H$ can
be explicitly determined by
\begin{equation}\label{eq:H}
H(\z)= \tanh (\frac{\sqrt{2}}{2} \z),
\end{equation}
and moreover the following estimates hold
\begin{equation}
\label{asyw}
  \left\{
\begin{array}{ll}
    H(\z)  - 1= - A_0 e^{- \sqrt{2} |\z|} + O(e^{-(2 \sqrt{2}) |\z|}) \ \
    \mbox{for} \ \z \to + \infty; & \\
        H(\z)  + 1=  A_0 e^{- \sqrt{2} |\z|} + O(e^{-(2 \sqrt{2}) |\z|})
        \ \ \mbox{for} \ \z \to - \infty; & \\
    H'(\z) =  \sqrt{2}  A_0 e^{- \sqrt{2} |\z|} + O(e^{-(2 \sqrt{2}) |\z|})
    \ \ \mbox{for} \ |\z| \to + \infty,  &  \\
\end{array}
\right.
\end{equation}
where $A_0>0$ is a fixed constant. We have the following well-known
result (we refer to Lemma 4.1 in \cite{mu} for the proof).

\begin{lem}
\label{linear} Consider the following eigenvalue problem
\begin{equation}
  \phi^{''} + (1-3H^2) \phi =\lambda \phi, \qquad \qquad
 \phi \in H^1 (\R).
\end{equation}
Then, letting $(\l_j)_{j}$ be the eigenvalues arranged in non
increasing order (counted with multiplicity) and $(\phi_j)_j$ be the
corresponding eigenfunctions, one has that
\begin{equation}
\lambda_1=0, \quad \phi_1= c H^{'}; \qquad \qquad \lambda_2 <0.
\end{equation}
As a consequence (by Fredholm's alternative), given any function
$\psi \in L^2(\R)$ satisfying $\int_{\R}  \psi H' =0$, the following
problem has a unique solution $\phi$
\begin{equation}
\phi^{''} + (1-3H^2) \phi =  \psi \ \hbox{ in }  \R, \qquad \quad
\int_{\R} H' \phi =0.
\end{equation}
Furthermore, there exists a positive constant $C$ such that
$\|\phi\|_{H^1(\R)} \leq C \|\psi\|_{L^2(\R)}$.
\end{lem}
We collect next some useful formulas: first of all we notice that
\begin{equation}\label{eq:H'form}
    H' = \frac{1}{\sqrt{2}} (1-H^2)\qquad \quad \hbox{ and } \qquad
    \quad H''=-\sqrt2 H H'.
\end{equation}
Moreover, setting
\begin{equation}\label{eq:L0}
    \mathcal{L}_0 u = u'' + (1-3H^2) u,
\end{equation}
we have that
\begin{equation}\label{eq:L0HH'}
    \mathcal{L}_0 (H H') = - 3 \sqrt{2} H (H')^2.
\end{equation}

\subsection{Geometric background}\label{ss:gb}
In this subsection we expand the coefficients of the metric in
local normal coordinates. We then derive as a consequence an expansion for the
Laplace-Beltrami operator. First of all, it is convenient
 to scale by $\frac 1 \e$ the coordinates in equation \eqref{eq:fife2} to obtain
\begin{equation}\label{eq:fife2e}
    \left\{
      \begin{array}{ll}
        \D u = (u-a(\e x))(u^2-1) & \hbox{ in } \O_\e, \\
        \frac{\pa u}{\pa \nu} = 0 & \hbox{ on } \pa \O_\e,
      \end{array}
    \right.
\end{equation}
where we have set $\O_\e = \frac 1 \e \O$. Following the same
notation we also set $K_\e = \frac 1 \e K$, and for $\g \in (0,1)$
we define
$$
  S_\e = \left\{ x \in \O_\e \; : \; dist(x,K_\e) < \e^{-\g}
\right\}; \qquad \qquad I_\e = [- \e^{-\g}, \e^{-\g}].
$$
We parameterize elements $x \in S_\e$ using their closest point $y$
in $K_\e$ and their distance $\z$ (with sign, positive in the
dilation of $\O_+$). Precisely, we choose a system of coordinates
$\ov{y}$ on $K$, and denote by ${\bf n}(\ov{y})$ the (unique) unit
normal vector to $K$ (at the point with coordinates $\ov{y}$)
pointing towards $\O_-$. Choosing also coordinates $y$ on $K_\e$
such that $\ov{y} = \e y$, we define the diffeomorphism $\G_\e :
K_\e \times I_\e \to S_\e$ by
\begin{equation}\label{eq:Gammae}
\G_\e (y,\z) = y + \z {\bf n} (\e y).
\end{equation}
We let the upper-case indices $A, B, C, \dots$ run from $1$ to
$n$, and the lower-case indices $ i, j, l, \dots$ run from $1$ to
$n-1$. Using some local coordinates $(y_i)_{i=1, \dots, n-1}$ on
$K_\e$, and letting $\var_\e$ be the corresponding immersion into
$\R^n$, we have
$$
  \begin{cases}
    \frac{\partial \G_\e}{\partial y_i} (y,\z) = \frac{\partial
    \var_\e}{\partial y_i} (y) + \e \z \frac{\partial {\bf n}}{\partial
    y_i} (\e y) = \frac{\partial \var_\e}{\partial y_i} (y) + \e \z
    \kappa^j_i (\e y) \frac{\partial \var_\e}{\partial y_j} (y), &
    \text{ for } i = 1, \dots, n-1; \\
    \frac{\partial \G_\e}{\partial \z} (y,\z) = {\bf n} (\e y). &
  \end{cases}
$$
where $(\kappa^j_i)$ are the coefficients of the mean-curvature
operator on $K$. Let also $(\ov{g}_{ij})_{ij}$ be the coefficients
of the metric on $K_\e$ in the above coordinates $y$. Then, letting
$g = g_\e$ denote the metric on $\O_\e$ induced by $\R^n$, we have
\begin{equation}\label{eq:gAB}
  g_{AB} = \left( \frac{\partial \G_\e}{\partial x_A},
  \frac{\partial \G_\e}{\partial x_B} \right) = \begin{pmatrix}
    (g_{ij}) & 0 \\ 0 & 1 \
  \end{pmatrix},
\end{equation}
where
\begin{eqnarray*}
  g_{ij} & = & \left( \frac{\partial \var_\e}{\partial y_i} (y) + \e \z
\kappa^k_i (\e y) \frac{\partial \var_\e}{\partial x_k} (y),
\frac{\partial \var_\e}{\partial y_j} (y) + \e \z \kappa^l_j (\e
y) \frac{\partial \var_\e}{\partial x_l} (y) \right) \\ & = &
\ov{g}_{ij} + \e \z \left( \kappa^k_i \ov{g}_{kj} + \kappa^l_j
\ov{g}_{il} \right) + \e^2 \z^2 \kappa^k_i \kappa^l_j \ov{g}_{kl}.
\end{eqnarray*}
Note that also the inverse matrix $\{g^{AB}\}$ decomposes as
$$
  g^{AB} = \begin{pmatrix}
    (g^{ij}) & 0 \\
    0 & 1 \
  \end{pmatrix}.
$$
From the above decomposition of $g_{AB}$ (and $g^{AB}$) and for $u$
defined on $S_\e$, one has
\begin{eqnarray}\label{eq:defLcoord} \nonumber
  \D_g u & = & g^{AB} u_{AB} + \frac{1}{\sqrt{\det g}}
\partial_A \left( g^{AB} \sqrt{\det g} \right) u_{B} \\ & = &
u_{\z \z} + g^{ij} u_{ij} +\frac{1}{\sqrt{\det g}}
\partial_\z \left( \sqrt{\det g} \right) u_\z +
\frac{1}{\sqrt{\det g}} \partial_i \left( g^{ij} \sqrt{\det g}
\right) u_j.
\end{eqnarray}
We have, formally
\begin{equation}\label{eq:detg}
  \det g = det (\ov{g}^{-1} g) \det \ov{g} = (\det \ov{g}) \left( 1
+ \e \z \tr (\ov{g}^{-1} \a) \right) + o(\e),
\end{equation}
where
$$
\a_{ij} = \kappa^k_i \ov{g}_{kj} + \kappa^l_j \ov{g}_{il}.
$$
There holds
$$
  (\ov{g}^{-1} \a)_{is} = \ov{g}^{sj} \a_{ij} = \ov{g}^{sj} \left(
  \kappa^k_i \ov{g}_{kj} + \kappa^l_j \ov{g}_{il} \right),
$$
and hence
\begin{equation}\label{eq:trga}
  \tr (\ov{g}^{-1} \a) = \ov{g}^{ij} \left( \kappa^k_i \ov{g}_{kj} +
  \kappa^l_j \ov{g}_{il} \right) = 2 \ov{g}^{ij} \kappa^k_i \ov{g}_{kj} = 2
  \kappa^i_i.
\end{equation}
We recall that the quantity $\kappa^i_i$ represents the mean
curvature of $K$, and in particular it is independent of the choice
of coordinates.

\

\no We note that the metric $g_{AB}$ can be expressed in function of
the metric $\ov{g}_{ij}$, the operator $\kappa^i_j$, and the
variable $\e \z$. Hence, fixing an integer $k$ and using a Taylor
expansion, we can write
\begin{equation}\label{eq:gdng}
  \frac{1}{\sqrt{\det g}} \partial_\z \sqrt{\det g} = \sum_{\ell=1}^k
  \e^\ell \z^{\ell-1} \tilde{G}_\ell(\e y) + \tilde{G} (\e y, \z),
\end{equation}
where $\tilde{G}_\ell : K\to \R$ are smooth functions, and
$\tilde{G}$ satisfies
\begin{equation}\label{eq:tG}
  \left\| \tilde{G}(\cdot, \z) \right\|'_{C^m(K)}
  \leq C_{k,m} |\z|^k \e^{k+1}, \qquad \qquad \z \in I_\e,
\end{equation}
where $C_{k,m}$ is a constant depending only on $K$, $k$,
and $m$. Again (and in the following), when we write $\| \cdot
\|'$ we keep the variable $\z$ fixed. In particular, from the
above computations it follows that
\begin{equation}\label{eq:tildeG1}
    \tilde{G}_1(\e y) = \kappa(\e y) := \kappa^i_i(\e y).
\end{equation}

\

\noindent We need now a similar expansion for the operator $\D_g$:
fixing the variable $\z \in I_\e$, the metric $g (y, \z) = g_\e (y,
\z)$ induces a metric $\hat{g}_{\e,\z}$ on $K$ in the following way.
Consider the homothety $i_\e : K\to K_\e$. We define
$\hat{g}_{\e,\z}$ to be
$$
\hat{g}_{\e, \z} = \e^2 i_\e^* g_\e (\cdot, \z),
$$
where $i_\e^*$ denotes the pull-back operator. Basically, we are
freezing the variable $\z$ and letting $y$ vary. Fixing an integer
$k$, for any smooth function $v : K\to \R$ we have the expansion
below, which follows from \eqref{eq:defLcoord}, reasoning as for
\eqref{eq:gdng}
\begin{equation}\label{eq:Li}
   \D_{\hat{g}_{\e, \z}} v =  \sum_{\ell=0}^k (\e \z)^\ell
  L_\ell v + \tilde{L}_{\e,k+1} v = \D_K v + \sum_{\ell=1}^k (\e \z)^\ell
  L_\ell v + \e^{k+1} \tilde{L}_{\e,k+1} v.
\end{equation}
Here $\{L_i\}_i, \tilde{L}_{\e,k+1}$ are linear second-order
differential operators acting on $y$ and satisfying
\begin{equation}\label{eq:Lie}
  \left\| L_i v \right\|'_{C^m(K)} \leq C_m
  \| v \|'_{C^{m+2}(K)}; \qquad \qquad
  \left\| \tilde{L}_{\e,k+1} v \right\|'_{C^m(K)} \leq C_m
  |\z|^{k+1} \| v \|'_{C^{m+2}(K)}
\end{equation}
for all smooth $v$, where $C_m$ is a constant depending only on
$K$, $k$, and $m$.

\

\noindent Consider now a function $u : S_\e \to \R$ of
the form
\begin{equation}\label{eq:utu}
  u(y, \z) = \tilde{u} (\e y, \z), \qquad \qquad \qquad y
  \in K_\e,\; \z \in I_\e.
\end{equation}
Then, scaling in the
variable $y$, we have
$$
   \D_{g_\e} u (y, \z) =  \tilde{u}_{\z \z} (\e y, \z) +
  \frac{1}{\sqrt{\det g}} \partial_\z \left( \sqrt{\det g} \right)
  \tilde{u}_\z (\e y, \z) + \e^2 \D_{\hat{g}_{\e, \z}} \tilde{u}
  (\e y, \z).
$$
Using the expansions \eqref{eq:gdng}, \eqref{eq:Li} together with
\eqref{eq:tildeG1}, the last equation becomes
\begin{eqnarray}\label{eq:dgu2}
  \D_{g_\e} u (y, \z) & = & \tilde{u}_{\z \z} (\e y, \z) +
  \left( \e \kappa + \sum_{\ell=2}^k \e^\ell \z^{\ell-1} \tilde{G}_\ell \right)
  \tilde{u}_\z (\e y, \z) + \tilde{G} (\e y, \z)
  \tilde{u}_\z (\e y, \z) \nonumber \\ & +& \e^2 \D_K \tilde{u}  + \sum_{\ell=1}^k
  \e^{2+\ell} \z^\ell L_\ell \tilde{u} (\e y, \z) + \e^{k+3} \tilde{L}_{\e,k+1} \tilde{u} (\e y, \z).
\end{eqnarray}

\subsection{Spectral analysis}\label{ss:sa}

We define the scaled Euler functional $\ov{J}_\e(u)$ in $\O_\e$
by
\begin{equation}\label{eq:Euler}
\ov{J}_\e(u)=\frac12\int_{\O_\e} |\n u|^2+\int_{\O_\e} F(\e
x,u)dx,\quad \hbox{ with}\quad F(x,u):=\int_{-1}^u
(s-a(x))(s^2-1)\,ds.
\end{equation}
We set for brevity
\begin{equation}\label{eq:bb}
b(y) := \pa_{\bf{n}} a(y,0)
\end{equation}
and we notice that by our choice of ${\bf n}$, we have $b> 0$ on
$K$. Now, we let $\var_j$ and $\l_j$ be the eigenfunctions and the
eigenvalues (with weight $b$) of
$$
  - \D_{K} \varphi_j = \l_j b(\ov{y}) \varphi_j.
$$
The $\l_j$'s can be obtained for example using the  Rayleigh
quotient: precisely if $M_j$ denotes the family of $j$-dimensional
subspaces of $H^1(K)$, then one has
\begin{equation}\label{eq:ljcourant}
 \l_j = \inf_{M \in M_j} \sup_{\varphi \in M, \varphi \neq 0}
 \frac{\int_{K} |\n_{K}\varphi|^2}{\int_{K} b(\overline{y}) \varphi^2}
 =\sup_{M \in M_{j-1}}\inf_{\varphi \perp M, \varphi \neq 0}
 \frac{\int_{K} |\n_{K}\varphi|^2}{\int_{K} b(\overline{y}) \varphi^2},
\end{equation}
where $\perp$ denotes the orthogonality with respect to the $L^2$
scalar
 product with weight $b$. We can estimate the $\l_j$ using a standard Weyl's
 asymptotic formula (\cite{chavel}), one has
\begin{equation}\label{eq:weyllj}
  \l_j \simeq C_{K,b} j^{\frac{2}{n-1}} \qquad \quad \hbox{ as } j \to
  + \infty,
\end{equation}
for some constant $C_{K,b}$ depending only on $K$ and $b$.

\section{Asymptotic analysis}\label{s:aa}

This section is devoted to the construction of approximate solutions
to \eqref{eq:fife2e}, and of approximate eigenfunctions (and
eigenvalues) in the $\z$ component (see the coordinates introduced
in \eqref{eq:Gammae}) of the relative linearized equation. Then we
characterize, via Fourier analysis, the profile of resonant
eigenfunctions in both the variables $y$ and $\z$.

\subsection{Approximate solutions and eigenfunctions}\label{ss:appsol} In this section,
given any integer $k\ge 1$, we construct an approximate solution
$u_{k,\e}$ to problem \eqref{eq:fife2e},  which solves the equation
up to an error of order $\e^{k+1}$. Using the above parametrization
$(y,\z)$ in $S_\e$, we make the following ansatz
\begin{equation}\label{eq:approx}
 u_{k,\e}(y,\z)=H(\z-\Phi(\e y))+\sum_{i=1}^k \e^i h_i(\e y,\z-\Phi(\e
 y)) \qquad \quad \hbox{ in } S_\e,
\end{equation}
where $H$ is the hetheroclinic solution of \eqref{ode} and where
$\Phi(\e y)=\Phi_0(\e y)+\sum_{i=1}^{k-1}\e^i\Phi_i(\e y)$ for some
smooth functions $\Phi_j$  defined on $K$. The corrections $(h_i)_i$
and $(\Phi_i)_i$ are to be constructed recursively in the index $i$,
depending on the Taylor expansion of $a$ and the geometry of $K$.
Since all the $h_i$'s will turn out to have an exponential decay in
$\z$, $u_{k,\e}$ can be easily extended (via some cutoff functions)
to an approximate solution in the whole $\O_\e$, see
\eqref{eq:hatuke} below.

We first determine $h_1$ by solving the equation up to an error of
order $\e^2$. To this aim, we expand the function $a$ in powers of
$\e$ as (notice that $a(\e y, 0) \equiv 0$)
\begin{equation}\label{eq:expa}
    a(\e y,\e \z) = \e b(\e y)\z + \sum_{\ell = 2}^k (\e \z)^\ell b_\ell(\e y)
    + \tilde{b}(y,z),
\end{equation}
where $\tilde{b}(y,z)$ is smooth and satisfies
$$
  |\tilde{b}(y,z)| \leq C_k |\e z|^{k+1}.
$$
Using the above expansion of the metric coefficients and the
Laplace-Beltrami operator, see in particular \eqref{eq:dgu2},
setting $s = \z - \Phi$, we obtain that the term (formally) of order
$\e$ in the equation is identically zero if and only if the
correction $h_1$ satisfies
\begin{equation}\label{eq:h1}
 \mathcal{L}_{0}h_1:=(h_1)_{ss}+(1-3H(s)^2)h_1=-\kappa(\e y) H'(s)+(s+\Phi_0) b(\e y) (1-H^2(s)).
\end{equation}
By the asymptotics in \eqref{asyw}, the right-hand side is of class
$L^2$ in $\R$ and, by Lemma \ref{linear}, \eqref{eq:h1} is solvable
provided the latter is orthogonal in $L^2$ to the function $H'(s)$.
Since $H'(s)$ is even in $s$ and since $b(\e y) > 0$, this is
possible choosing $\Phi_0(\e y)$ so that
\begin{equation}\label{eq:bPhi0}
    b(\e y) \Phi_0(\e y) =\kappa(\e y) \frac{\int_{-\infty}^{+\infty}(H')^2(s)
    ds}{\int_{-\infty}^{+\infty}H'(s)(1-H^2(s)) ds}=\frac{\sqrt2}{3}\k(\e y).
\end{equation}
Moreover one can prove, using standard ODE estimates, that $h_1$ has
the following (regularity properties and) decay at infinity
\begin{equation}\label{eq:dech1}
    |\pa_s^l \pa_y^m h_1(\e y,s)| \leq C_m \e^m (1 + |s|) e^{-\sqrt{2}
    |s|}, \qquad \quad l = 0, 1, 2, \quad m = 0, 1, 2, \dots,
\end{equation}
where $C_m$ depends only on $m$, $a$ and $K$.

\

\noindent To obtain the other corrections $(\Phi_i)_i$ and $(h_i)_i$
one can proceed by induction, assuming that $N \geq 2$, that
$\Phi_0, \dots, \Phi_{N-2}$ and $h_1, \dots, h_{N-1}$ have been
determined, and that $(h_i)_{i\le N-1}$ satisfy
\begin{equation}\label{eq:dechi}
    |\pa_s^l \pa_y^m h_i(\e y,s)| \leq C_m \e^m (1 + |s|^{d_i}) e^{-\sqrt{2}
    |s|}, \qquad \quad i\le N-1\qquad l = 0, 1, 2, \quad m = 0, 1, 2, \dots,
\end{equation}
where $C_m$ depends only on $m$, $a$, $K$ and $d_i$ only on $i$.
When we expand the equation \eqref{eq:fife2e} for $u=u_{N,\e}$ in
power series of $\e$, the couple $(h_N,\Phi_{N-1})$ can be found
reasoning as for $(h_1,\Phi_0)$: indeed, considering the coefficient
of $\e^N$ in this expansion, one can easily see that $h_N$ satisfies
an equation on the form
\begin{equation}\label{eq:hN}
 \mathcal{L}_{0}h_N=\Phi_{N-1} b(\e y) (1-H^2(s))+
 \mathfrak{F}_N(s,\Phi_0,\dots,\Phi_{N-2},h_1,\dots,h_{N-1},b_1,\dots,b_N),
\end{equation}
where $\mathfrak{F}_N$ is a smooth function on its argument.
Reasoning as for $h_1$ this equation is solvable provided the
right-hand side is $L^2$-orthogonal to the function $H'(s)$. This is
indeed true choosing $\Phi_{N-1}$ so that
\begin{equation*}
    b(\e y) \Phi_{N-1}(\e y) =\frac{-\int_{-\infty}^{+\infty}H'(s)
    \mathfrak{F}_N(s,\cdot) ds}{\int_{-\infty}^{+\infty}H'(s)(1-H^2(s)) ds}.
\end{equation*}
Furthermore, one can show that $h_N$ satisfies regularity and decay
estimates as in \eqref{eq:dechi}. Reasoning as in Section 3 of
\cite{mm2} one can check that the above formal estimates can be made
rigorous, and that the exponential decay of the corrections yields
the following result.

\begin{pro}\label{p:as}
Given any integer $k\ge 1$ there exist a function
$u_{k,\e}:\,S_\e\to \R$ which solves equation \eqref{eq:fife2e} up
to an error of order $\e^{k+1}$. Precisely, setting
\begin{equation}\label{eq:mfse}
       \mathfrak{S}_\e(u)= \D u - (u-a(\e x))(u^2-1),
\end{equation}
there exist a polynomial $P_k(\z)$ such that
\begin{equation}\label{eq:deceqas}
  |\mathfrak{S}_\e(u_{k,\e}(\e y, \z))| \leq \e^{k+1} P_k(\z)
  e^{-\sqrt{2} |\z|} \qquad \quad \hbox{ in } S_\e.
\end{equation}
Moreover, the following estimate holds
\begin{equation}\label{eq:decuke}
    |\pa_s^l \pa_y^m u_{k,\e}(\e y,\z)| \leq C_m \e^m P_k(\z) e^{-\sqrt{2}
    |\z|}, \qquad \quad l = 0, 1, 2, \quad m = 0, 1, 2, \dots,
\end{equation}
where $C_m$ is a constant depending only on $m$, $a$ and $K$.
\end{pro}

\

\no We will look at solutions $u$ of \eqref{eq:fife2e} as small
corrections of $u_{k,\e}$ (suitably extended to $\O_\e$ via some
cutoffs in $\z$, see \eqref{eq:hatuke} below), namely of the form
 $$
 u=u_{k,\e}+w,
 $$
for $w$ small in a sense to be specified later. The equation
$\mathfrak{S}_\e(u)=0$ is then equivalent to
$$
L_\e(w)+\mathfrak{S}_\e(u_{k,\e})+\mathcal{N}_\e(w)=0,
$$
where $L_\e$ is nothing but the linearized operator at the
approximate solution $u_{k,\e}$
 \begin{equation}\label{eq:linop}
  L_\e w:=\D_{g_\e} w+(1-3u_{k,\e}^2)w+2a(\e x)u_{k,\e} w ,
 \end{equation}
 and where $\mathcal{N}_\e$ is the remainder given by nonlinear
 terms in $\mathfrak{S}_\e$, namely
\begin{equation}\label{eq:nlinop}
   \mathcal{N}_\e(w):=-(3 u_{k,\e}-a(\e x))w^2-w^3.
\end{equation}
It is also convenient to define the following linear operator
\begin{equation}\label{mbble}
\mathbb{L}_\e w := w_{\z \z} + \frac{1}{\sqrt{\det g}} \partial_\z
\left( \sqrt{\det g} \right)
  w_\z + (1-3u_{k,\e}^2) w + 2 a(\e y, \z) u_{k,\e} w.
\end{equation}
In particular, using the expansion \eqref{eq:dgu2}, the operators
$L_\e$ and $\mathbb{L}_\e$ are related by the following formula:
setting $w(y,\z)=\tilde{w}(\e y,\z)$ one has
\begin{equation}\label{eq:expLe}
L_\e w = \mathbb{L}_\e \tilde{w} +\e^2\D_{K} \tilde{w} + \e^3
\hat{\mathfrak{L}}_{3,\e} \tilde{w},
\end{equation}
where $\hat{\mathfrak{L}}_{3,\e}$ consists of the last two terms in
\eqref{eq:dgu2} (replacing $u$ with $w$). Precisely,
$\hat{\mathfrak{L}}_{3,\e}$ is a linear differential operator of
second order acting on the variables $\ov{y}$, which for every
integer $m$ satisfies
\begin{equation}\label{eq:Lja}
  \left\| \hat{\mathfrak{L}}_{3,\e} v \right\|'_{C^m(K)} \leq
  \hat{P}(\z) \| v \|'_{C^{m+2}(K)}.
\end{equation}
Here $\hat{P}(\z)$ is a polynomial in $\z$ with fixed degree, and
coefficients depending only on $m$.

\

We want next to derive some formal estimates on the following
eigenvalue problem
\begin{equation}\label{eq:linn}
   \mathbb{L}_\e  v = \mu v \qquad
  \hbox{ in } I_\e,
\end{equation}
with zero Dirichlet boundary conditions. It follows from Lemma
\ref{linear} that the eigenvalues either stay bounded away from
zero, or converge to zero as $\e \to 0$: we are interested in the
latter case. We argue heuristically expanding \eqref{eq:linn} at
first order in $\e$. In the limit $\e \to 0$, we have $\mu = 0$ with
corresponding eigenfunction $H'$, therefore it is natural to look
for approximate eigenfunctions of the form
\begin{equation}\label{eq:Psi}
\Psi = H' + \e H_1,
\end{equation}
and eigenvalues $\mu = \e \ov{\mu} + o(\e)$. We impose that $H_1$ is
orthogonal to $H'$ in $L^2(\R)$. Therefore the approximate
eigenvalue equation formally becomes
\begin{equation}\label{eq:H1}
  \mathcal{L}_0 H_1 = - 2 b(\e y) (s + \Phi_0) H H' - H'' \k(\e y) + 6
  H H' h_1 + \ov{\mu} H' + o(1).
\end{equation}
As for \eqref{eq:h1}, solvability is guaranteed provided the
right-hand side is orthogonal in $L^2$ to $H'$. Using the oddness of
$H$, formulas \eqref{eq:H'form}, \eqref{eq:L0HH'}, \eqref{eq:h1} and
the self-adjointness of $\mathcal{L}_0$ we find that orthogonality
is equivalent to
$$
  \ov{\mu} = 4 b(\e y) \frac{\int_\R s H(s) (H'(s))^2 ds}{\int_\R (H'(s))^2
  ds}  = \sqrt{2} b(\e y).
$$
With this choice of $\ov{\mu}$, the function $H_1$ is defined as the
unique solution of
$$
  \mathcal{L}_0 H_1 = - 2 b(\e y) (s + \Phi_0) H H' - H'' \k(\e y) + 6
  H H' h_1 + \sqrt{2} b(\e y) H'.
$$
From the exponential decay of $H'$ and $h_1$ (see \eqref{eq:dech1})
we deduce that $H_1$ satisfies estimates similar to
\eqref{eq:dechi}.

Using this fact and the rigorous expansions in \eqref{eq:gdng},
\eqref{eq:expa}, we then derive the following estimate
\begin{equation}\label{eq:approxeig}
    \mathbb{L}_\e(\Psi) = \e \ov{\mu} H' + \e^2 R_\e(\e y,
    \z),
\end{equation}
where the error term $R_\e$ satisfies
\begin{equation}\label{eq:Re}
    |R_\e(\e y, \z)| \leq P(\z)e^{-\sqrt{2}|\z|}
\end{equation}
for some polynomial $P(\z)$. Also, from the regularity and the decay
of $H_1$ we have
\begin{equation}\label{eq:decPsi}
    |\pa_s^l \pa_y^m \Psi(\e y,s)| \leq P(\z) \e^{m+1} e^{-\sqrt{2}
    |\z|}, \qquad \quad l = 0, 1, 2, \quad m = 0, 1, 2, \dots.
\end{equation}

\subsection{Characterization of resonant
eigenfunctions}\label{ss:re}

We characterize next the eigenfunctions of $L_\e$, see
\eqref{eq:linop}, corresponding to small eigenvalues. Let us recall
first the definition of $(\varphi_j)_j$ and $(\l_j)_j$ in Subsection
\ref{ss:sa}.

\begin{lem}\label{l:eigenlinop} Let $\l_\e=O(\e^2)$ be an eigenvalue of the
linearized operator $L_\e$ in $S_\e$ with eigenfunction $\phi$ and
weight $b$, namely
$$
L_\e\phi=\l_\e b \phi \qquad \quad \hbox{ in } S_\e
$$
(and with zero Dirichlet boundary conditions). Let us write the
eigenfunction $\phi$ as
$$
\phi(y,\z)=\,\var(\e y)\Psi(\e y,\z)+\phi^\perp(y,\z),
$$
with $\Psi$ defined in \eqref{eq:Psi} and with $\phi^\perp$
satisfying the following orthogonality condition (we are freezing
the $y$ variables in the volume element)
\begin{equation}\label{eq:orth}
\int_{I_\e} \Psi(\e y,\z) \phi^\perp(y,\z)\, dV_{g_\e}(\z) = 0\qquad
\hbox{for every}\quad y\in K_\e.
\end{equation}
Then  one has $\|\phi^\perp\|_{H^1(S_\e)}=o(\e)
\|\phi\|_{H^1(S_\e)}$ as $\e$ tends to zero.
\end{lem}

\begin{pf}
We notice first that, since $\Psi = H' + o(1)$ in $H^1(\R)$, the
operator $\mathcal{L}_0$ is negative definite on $\phi^\perp$ by
Lemma \ref{linear}. Therefore, using the estimates on the metric
$g_\e$ in Section \ref{s:not}, we find easily that there exist a
constant $C>0$ such that
\begin{equation}\label{eq:posphip}
\int_{S_\e}\phi^\perp(y,\z) L_\e \phi^\perp(y,\z)\,dV_{g_\e}(y,\z)
\le - \frac 1 C \|\phi^\perp\|^2_{H^1(S_\e)}.
\end{equation}
Let us write the eigenvalue equation $L_\e\phi=\l_\e b \phi$ as
$$
L_\e\phi^\perp=-L_\e(\var\Psi)+\l_\e b \phi^\perp + \l_\e b
\var\Psi.
$$
Multiplying by $\phi^\perp$, integrating over $S_\e$ and using
\eqref{eq:orth} we obtain
\begin{eqnarray*}
 \int_{S_\e}\phi^\perp(y,\z) L_\e \phi^\perp(y,\z)\,dV_{g_\e}&=&
- \int_{S_\e}\phi^\perp(y,\z) L_\e (\var(\e y)\Psi(\e
y,\z))\,dV_{g_\e} +\l_\e\int_{S_\e} b(\e y) \phi^\perp(y,\z)^2
\,dV_{g_\e}.
\end{eqnarray*}
By \eqref{eq:posphip} (and the
smallness of $\l_\e$) it then follows
\begin{equation}\label{eq:fin}
    \|\phi^\perp\|^2_{H^1(S_\e)} \leq C \left| \int_{S_\e}\phi^\perp(y,\z)
    L_\e (\var(\e y)\Psi(\e y,\z))\,dV_{g_\e} \right|.
\end{equation}
Now by \eqref{eq:expLe} and \eqref{eq:approxeig} we can write
$L_\e(\var\Psi)$ as
\begin{eqnarray}\label{eq:LembbLe} \nonumber
  L_\e(\var\Psi) & = & \mathbb{L}_\e(\var\Psi)+\e^2
  \D_{K}(\var\Psi)+ \e^3\hat{\mathfrak{L}}_{3,\e}(\var\Psi) \\
  & = & \var\left( \e\sqrt{2}b(\e y)H'+\e^2 R_\e(\e y,\z)\right)
+\e^2 \D_{K}(\var\Psi)+ \e^3\hat{\mathfrak{L}}_{3,\e}(\var\Psi).
\end{eqnarray}
Then, again by \eqref{eq:orth}, we have
 \begin{eqnarray}\label{eq:phiperpLePp}
\nonumber
  \left| \int_{S_\e}\phi^\perp(y,\z) L_\e (\var(\e y)\Psi(\e y,\z)) \right| &\le&
 \e^2 \left|\int_{S_\e}  \phi^\perp(y,\z)\D_K(\var(\e y)\Psi(\e y,\z))\right|\\
 &+&\left|\int_{S_\e} \left(\e^2 \tilde{R}_\e(\e y,\z)\var(\e y)
+ \e^3\hat{\mathfrak{L}}_{3,\e}(\var(\e y)\Psi(\e
y,\z))\right)\phi^\perp\right|,
 \end{eqnarray}
where $\tilde{R}_\e$ is as in \eqref{eq:Re}. We first estimate the
second term: since $\hat{\mathfrak{L}}_{3,\e}$ is a second order
operator in $\ov{y}$ satisfying the bound \eqref{eq:Lja}, by an
integration by parts we find
$$
  \e^3 \left|\int_{S_\e} \hat{\mathfrak{L}}_{3,\e}(\var(\e y)\Psi(\e
  y,\z)) \phi^\perp\,dV_{g_\e}(y,\z) \right| \leq C \e^2 \int_{K_\e \times
  I_\e} \hat{P}(\z) |\n_{\ov{y}} (\var(\e y)\Psi(\e
  y,\z))| \, |\n_y \phi^\perp| dV_{g_\e}(y,\z).
$$
Therefore, using the H\"older inequality, the change of variables
$\ov{y} = \e y$, \eqref{eq:decPsi} and the estimate on
$\tilde{R}_\e$ we find
\begin{equation}\label{eq:esterror}
\left|\int_{S_\e} \left(\e^2 \tilde{R}_\e(\e y,\z)\var(\e y) +
\e^3\hat{\mathfrak{L}}_{3,\e}(\var(\e y)\Psi(\e y,\z))\right)
\phi^\perp\,dV_{g_\e}(y,\z) \right|\le C
\frac{\e^2}{\e^{\frac{n-1}{2}}}
 \|\phi^\perp\|_{H^1({S_\e})} \|\var\|_{H^1(K)},
\end{equation}
for some positive constant $C$. It remains to estimate the first
term in \eqref{eq:phiperpLePp}. To this aim we decompose $\var$ as
(see the above notation)
\begin{equation}\label{eq:deccphi}
    \var(\e y)=\sum_j \a_j\varphi_j(\e y),
\end{equation}
for some real numbers $\a_j$. One can write
\begin{eqnarray*}
   \int_{S_\e}  \phi^\perp(y,\z)\D_K(\var(\e y)\Psi(\e y,\z)) \,dV_{g_\e}(y,\z) &=& \int_{S_\e}
   \phi^\perp(y,\z)\Psi(\e y,\z)\D_K\var(\e y)\,dV_{g_\e}(y,\z)\\
   &+&\int_{S_\e}  \phi^\perp(y,\z)\var(\e y)\D_K \Psi(\e y,\z)\,dV_{g_\e}(y,\z)\\
   &+&2\int_{S_\e}  \phi^\perp(y,\z)\n_K\Psi(\e y,\z)\cdot\n_K\var(\e y)\,dV_{g_\e}(y,\z).
\end{eqnarray*}
The first term vanishes by \eqref{eq:orth}. Hence, using the
H\"older inequality, \eqref{eq:decPsi} and reasoning as for
\eqref{eq:esterror} we obtain
\begin{eqnarray*}
\left| \int_{S_\e} \phi^\perp(y,\z)\D_K(\var(\e y)
\Psi(y,\z))\,dV_{g_\e}(y,\z)
 \right| \leq C \frac{1}{\e^{\frac{n-1}{2}}}
  \|\phi^\perp\|_{H^1({S_\e})} \|\var\|_{H^1(K)},
\end{eqnarray*}
which by estimates \eqref{eq:phiperpLePp} and \eqref{eq:esterror}
implies
\begin{equation}\label{eq:phiperpLe}
    \left| \int_{S_\e}\phi^\perp(y,\z) L_\e (\var(y)\Psi(y,\z)) \,
    dV_{g_\e}(y,\z)\right| \leq C \frac{\e^2}{\e^{\frac{n-1}{2}}}  \|\phi^\perp\|_{H^1({S_\e})}
    \|\var\|_{H^1(K)}.
\end{equation}
Using then \eqref{eq:fin} and the latter equation together with the
Weyl's asymptotic formula (see Subsection \ref{ss:sa}), we then
obtain
\begin{equation}\label{eq:estphiperp}
    \|\phi^\perp\|_{H^1(S_\e)}^2 \leq \frac{C}{\e^{n-1}} \left(
    \sum_j\a_j^2(\e^4 + \e^4 j^{\frac{2}{n-1}}) \right).
\end{equation}
Now, we rewrite the eigenvalue equation as
\begin{eqnarray*}
\left(\e\sqrt2\,b\,\var+\e^2\D_{K}  \var \right)\Psi &=& \l_\e b
\phi^\perp+\l_\e b \var\Psi -\e^2\var\D_{K} \Psi -
2\e^2\n_K\var\cdot\n_K\Psi\\
&-&\e^2 \tilde{R}_\e \var -\e^{3}\hat{\mathfrak{L}}_{3,\e}(\var
\Psi)- L_\e\phi^\perp.
\end{eqnarray*}
We use again the above decomposition $\varphi(\e y) = \sum_j \a_j
\varphi_j(\e z)$, we define the integer $j_\e$ (depending on $\e$)
to be the first $j$ such that $\e^2 \l_j > \e^{\frac 12}$. We
multiply this time the last equation by $\sum_{j \geq j_\e} \a_j
\varphi_j \Psi$ and we integrate over $S_\e$. By  the
self-adjointness of $L_\e$, \eqref{eq:decPsi}, \eqref{eq:orth} and a
similar argument as for \eqref{eq:esterror}, incorporating the term
involving $\hat{\mathfrak{L}}_{3,\e}$ into the left-hand side we
obtain that
\begin{eqnarray}\label{eq:fethi}\nonumber
  & & (1 + O(\e)) \frac{1}{\e^{n-1}} \sum_{j \geq j_\e} \e^2 \a_j^2  \l_j
  \leq C (\e^2 + |\l_\e|)
  \left( \frac{1}{\e^{n-1}} \sum_{j \geq j_\e} \a_j^2 \right)^{\frac 12}
  \left( \frac{1}{\e^{n-1}} \sum_{j} \a_j^2 \right)^{\frac 12}  \\
  & + & C \e \left( \frac{1}{\e^{n-1}} \sum_j \a_j^2 \right)^{\frac
  12} \left( \frac{1}{\e^{n-1}} \sum_{j \geq j_\e} \e^2  \l_j \a_j^2
  \right)^{\frac 12} + \left| \int_{S_\e} \phi^\perp L_\e
  (\sum_{j \geq j_\e} \a_j \varphi_j \Psi) dV_{g_\e} \right|.
\end{eqnarray}
The last term can be estimated reasoning as for
\eqref{eq:estphiperp}: since $\l_j \gg 1$ for $j \geq j_\e$, we find
that
\begin{eqnarray*}
  \left| \int_{S_\e} \phi^\perp L_\e (\sum_{j \geq j_\e} \a_j
  \varphi_j \Psi) dV_{g_\e} \right| \leq  C \e \|\phi^\perp\|_{H^1(S_\e)} \left(
  \frac{1}{\e^{n-1}} \sum_{j \geq j_\e} \e^2 \l_j \a_j^2 \right)^{\frac
  12}.
\end{eqnarray*}
From the above estimates and the fact that $\l_\e = O(\e^2)$ we
obtain
\begin{equation}\label{eq:ffff}
    \left( \frac{1}{\e^{n-1}} \sum_{j \geq j_\e} \e^2 \l_j \a_j^2
  \right)^{\frac 12} \leq C \e \left( \left( \frac{1}{\e^{n-1}} \sum_{j}
  \a_j^2 \right)^{\frac 12} + \|\phi^\perp\|_{L^2(S_\e)} \right).
\end{equation}
Then, writing $\var\Psi$ as
$$\var\Psi=\sum_{j < j_\e} \a_j\var_j\Psi+\sum_{j \ge
j_\e} \a_j\var_j\Psi,
$$
using \eqref{eq:estphiperp} and \eqref{eq:ffff} we find
\begin{eqnarray*}
  \|\phi^\perp\|_{H^1(S_\e)} & \leq & C (\e^2 + \e^{\frac 54})
  \|\var \Psi\|_{L^2(S_\e)} + C \e \left(
  \frac{1}{\e^{n-1}} \sum_{j \geq j_\e} \e^2 \l_j \a_j^2
  \right)^{\frac 12}  \\ & \leq & C \e^{\frac 54} \|\var \Psi\|_{L^2(S_\e)}
  + C \e^2 (\|\var \Psi\|_{L^2(S_\e)} + \|\phi^\perp\|_{H^1(S_\e)}).
\end{eqnarray*}
This implies $\|\phi^\perp\|_{H^1(S_\e)} \leq C \e^{\frac 54} \|\var
\Psi\|_{L^2(S_\e)}$, and noticing that $\|\phi\|_{L^2(S_\e)}^2 =
\|\var \Psi\|_{L^2(S_\e)}^2 + \|\phi^\perp\|_{L^2(S_\e)}^2$, we
achieve the desired estimate.
\end{pf}

\

\no Our next task is to estimate the derivatives of small
eigenvalues of the linearized operator $L_\e$ with respect to the
parameter $\e$. This will allow us to obtain invertibility of $L_\e$
for a suitable family of small $\e$. The prove of the main result
can be then obtained by a direct application of the contraction
mapping theorem. Using a result by T.Kato, see \cite{ka}, page 445,
which can be applied by the symmetry of $L_\e$ and elliptic
regularity results (these ensure that the eigenvalues of $L_\e$ are
{\em stable} and {\em semi-simple}, according to the definitions in
\cite{ka}), we have the following proposition.

\begin{pro}\label{p:epS}
For $k \in \N$, $k \geq 1$, let $u_{k,\e}$ be given by Proposition
\ref{p:as}, and let $L_\e$ be defined in \eqref{eq:linop}. Then the
eigenvalues $\l_\e$ of the problem
\begin{equation}\label{eq:epS}
  \left\{
    \begin{array}{ll}
      L_\e u = \l_\e b u & \hbox{ in } S_\e; \\
      u = 0 & \hbox{ on } \pa S_\e,
    \end{array}
  \right.
\end{equation}
are differentiable with respect to $\e$, and they satisfy the
following estimates
\begin{equation}\label{eq:estdle}
    T^1_{\l_\e,\e}  \leq \frac{\partial \l_\e}{\partial
    \e} \leq T^2_{\l_\e,\e}.
\end{equation}
Here we have set
$$
  T^1_{\l,\e} = \inf_{u \in H_\l, u \neq 0} \frac{\int_{S_\e} \left(-\frac 2 \e |\n_{g_\e}
    u|^2 - 6 u_{k,\e} v_{k,\e} u^2\right) dV_{g_\e}}{\int_{S_\e} b u^2 dV_{g_\e}};
$$
$$
  T^2_{\l,\e} = \sup_{u \in H_\l, u \neq 0} \frac{\int_{S_\e} \left(-\frac 2 \e  |\n_{g_\e}
    u|^2 - 6 u_{k,\e} v_{k,\e} u^2
    \right) dV_{g_\e}}{\int_{S_\e} b u^2 dV_{g_\e}},
$$
$v_{k,\e}=\frac{\pa u_{k,\e}}{\pa \e}$, while $H_\l$ stands for the
eigenspace for \eqref{eq:epS} corresponding to the eigenvalue $\l$.
\end{pro}

\noindent We next give a further characterization of some
eigenfunctions of $L_\e$, in addition to the ones in Lemma
\ref{l:eigenlinop}, concerning in particular the function $\varphi$.

\begin{lem}\label{l:charvarphi}
Suppose the assumptions of Lemma \ref{l:eigenlinop} hold true. Then,
normalizing $\phi$ by $\|\phi\|_{H^1(S_\e)} = 1$, decomposing $\var$
as in \eqref{eq:deccphi} and setting
\begin{equation}\label{eq:lje}
    \l_{j,\e} = \sqrt{2} \e - \e^2 \l_j,
\end{equation}
as $\e \to 0$ we have that
$$
  \frac{1}{\e^{n-1}}\sum_{|\l_{j,\e}| \geq \e^{\frac 54}} \a_j^2 = o(1);
  \qquad  \quad \frac{1}{\e^{n-1}} \sum_{|\l_{j,\e}| \geq \e^{\frac 54}} \l_{j,\e}
  \a_j^2 = o(\e).
$$
\end{lem}

\begin{pf}
We define the sets
$$
  A_{1,\e} = \left\{ j \in \N \; : \; \l_{j,\e} > \e^{\frac 54}
  \right\}; \qquad \qquad A_{2,\e} = \left\{ j \in \N \; : \;
  \l_{j,\e} < - \e^{\frac 54} \right\},
$$
and the functions
$$
  \ov{\varphi}_1(\e y) = \sum_{j \in A_{1,\e}} \a_j
  \varphi_j(\e y); \qquad \qquad \ov{\varphi}_2(\e y) = \sum_{j \in
  A_{2,\e}} \a_j \varphi_j(\e y);
$$
$$
   \phi_1 = \ov{\varphi}_1(\e y) \Psi(\e y,\z); \qquad \qquad \qquad
   \phi_2 = \ov{\varphi}_2(\e y) \Psi(\e y,\z).
$$
As one can easily see, from the estimates on $g_\e$ in Subsection
\ref{ss:gb} and from the decay of $\Psi$, see \eqref{eq:decPsi}, as
$\e \to 0$ there holds
\begin{equation}\label{eq:vvv}
    \|\var\Psi\|_{H^1(S_\e)}^2 = \frac{1+o(1)}{\e^{n-1}} \left( C_0 \int_K \var^2
  + C_1 \int_K |\n \var|^2 \right),
\end{equation}
where $C_0 = \int_\R (H')^2 + (H'')^2$ and $C_1 = \int_\R (H')^2 $.
Similar formulas hold true for $\phi_1$ and $\phi_2$, and hence
these two functions stay uniformly bounded in $H^1(S_\e)$ as $\e$
tends to zero.

We multiply next the equation in \eqref{eq:epS} by $\phi_l$, $l = 1,
2$: from the orthogonality of $\ov{\var}_1$ and $\ov{\var}_2$ on $K$
(with weight $b$), \eqref{eq:orth}, an integration by parts and the
above arguments we get
\begin{eqnarray*}
O(\e^2)\|\phi_l\|^2_{L^2(S_\e)}+ O(\e^3)\|\var\Psi\|_{L^2(S_\e)}
\|\phi_l\|_{L^2(S_\e)} =
  \int_{S_\e} \phi_l L_\e \phi dV_{g_\e} = \int_{S_\e} (\varphi \Psi +
  \phi^\perp) L_\e \phi_l.
\end{eqnarray*}
Using \eqref{eq:phiperpLe} (replacing $\var$ by $\ov{\var}_l$) we
deduce
$$
  \left| \int_{S_\e}\phi^\perp(y,\z) L_\e (\ov{\var}_l(y)\Psi(y,\z)) \,
    dV_{g_\e}(y,\z)\right| \leq C \frac{\e^2}{\e^{\frac{n-1}{2}}}
    \|\phi^\perp\|_{H^1({S_\e})} \|\ov{\var}_l\|_{H^1(K)},\qquad \hbox{ for } l=1,2.
$$
Also, from the expression of $L_\e$ (see \eqref{eq:LembbLe}), the subsequent estimates and some
straightforward computations one finds
\begin{eqnarray*}
  \int_{S_\e} \varphi \Psi L_\e \phi_l & = & - \frac{1 + o(1)}{\e^{n-1}}
  \sum_{j \in A_{l,\e}}
   \a_j^2 \l_{j,\e} + O(\e^2) \left( \frac{1}{\e^{n-1}}
   \sum_{j \in A_{l,\e}} \a_j^2 \right)^{\frac 12} \|\phi\|_{H^1(S_\e)}
   \\ & + & O(\e) \left( \frac{1}{\e^{n-1}} \sum_{j \in A_{l,\e}}
   \a_j^2 \e^2 \l_{j} \right)^{\frac 12} \|\phi\|_{H^1(S_\e)}.
\end{eqnarray*}
Then from the last three formulas,  Lemma \ref{l:eigenlinop} and the
normalization on $\phi$ we obtain
\begin{equation}\label{eq:fffin}
  \frac{1}{\e^{n-1}} \sum_{j \in A_{l,\e}} \a_j^2 \l_{j,\e} =
  O(\e^2) \left( \frac{1}{\e^{n-1}}
   \sum_{j \in A_{l,\e}} \a_j^2 \right)^{\frac 12}
   + O(\e) \left( \frac{1}{\e^{n-1}} \sum_{j \in A_{l,\e}}
   \a_j^2 \e^2 \l_{j} \right)^{\frac 12},
\end{equation}
for $l = 1, 2$.

\

\noindent We next further split the sets $A_{l,\e}$ as $A_{l,\e} =
\hat{A}_{l,\e} \cup \tilde{A}_{l,\e}$, where
$$
  \hat{A}_{l,\e} = \left\{ j \in A_{l,\e} \; : \; |\e^2 \l_j| < \e^{\frac 34}
  \right\}; \qquad \qquad \tilde{A}_{l,\e} = A_{l,\e} \setminus
  \hat{A}_{l,\e},
$$
so from \eqref{eq:fffin} we have clearly
\begin{eqnarray}\label{eq:ciaociao} \nonumber
  \frac{1}{\e^{n-1}} \sum_{j \in A_{l,\e}} \a_j^2 \l_{j,\e} & = &
  O(\e^2) \left( \frac{1}{\e^{n-1}}
   \sum_{j \in A_{l,\e}} \a_j^2 \right)^{\frac 12}
   + O(\e) \left( \frac{1}{\e^{n-1}} \sum_{j \in \hat{A}_{l,\e}}
   \a_j^2 \e^2 \l_{j} \right)^{\frac 12} \\ & + &
   O(\e) \left( \frac{1}{\e^{n-1}} \sum_{j \in \tilde{A}_{l,\e}}
   \a_j^2 \e^2 \l_{j} \right)^{\frac 12}.
\end{eqnarray}
Obviously, since $\l_{j,\e} = \sqrt{2} \e - \e^2 \l_j$, for $j \in
\tilde{A}_{l,\e}$ the ratio $\frac{\e^2 \l_j}{\l_{j,\e}}$ stays
uniformly bounded from above and below by positive constants.
Therefore, using the elementary inequality $|xy| \leq \d |x|^2 +
\frac{1}{\d} |y|^2$ with $\d$ small and fixed, we can absorb the
last term into the left-hand side of the latter formula, obtaining
an error of the form $O(\e^2/\d)$. Therefore, using also the
definition of $\hat{A}_{l,\e}$ we find
$$
  \frac{1}{\e^{n-1}} \sum_{j \in A_{l,\e}} \a_j^2 \l_{j,\e} =
  O(\e^{\frac{11}{8}}) \left( \frac{1}{\e^{n-1}}
   \sum_{j \in A_{l,\e}} \a_j^2 \right)^{\frac 12}+O(\e^2).
$$
By our normalization on $\phi$, see also the comments after
\eqref{eq:vvv}, the argument inside the last bracket is uniformly
bounded as $\e \to 0$, and hence we obtain the second assertion of
the lemma.

\

\no To obtain also the first one we notice that $|\l_{j,\e}| \geq
\e^{\frac 54}$ for $j \in A_{l,\e}$, so we find
$$
  \frac{1}{\e^{n-1}} \sum_{j \in A_{l,\e}} \a_j^2 = O(\e^{\frac{11}{8} -
  \frac{5}{4}}) = O(\e^{\frac{1}{8}}) = o(1)
$$
as $\e \to 0$. This concludes the proof.
\end{pf}

Now, using the above lemma, we can estimate the derivatives of small
eigenvalues of $L_\e$ with respect to $\e$. Precisely, we have the
following result.

\begin{lem}\label{l:dereigenv} Let $\l$ be as in Lemma \ref{l:eigenlinop}.
Then, for $\e$ sufficiently small $\l$ is differentiable with
respect to $\e$, and there exists a negative constant $c_{K,b}$,
depending only on $K$ and $b$, such that its derivative (which is
possibly a multi-valued function) satisfies
$$
  \left| \frac{\partial \l}{\partial \e} - c_{K,b} \right| = o(1)
  \qquad \hbox{ as } \e \to 0.
$$
\end{lem}

\begin{pf}
The proof is based on Lemma \ref{l:eigenlinop}, Proposition
\ref{p:epS} and Lemma \ref{l:charvarphi}. Since we want to apply
formula \eqref{eq:estdle} (in our previous notation) to the function
\[u= \phi = \left(\sum_{ |\l_{j,\e}| \geq \e^{\frac 54}} \a_j  \var_j(\e z)\right)
\Psi + \phi^\perp=\var \left(H'+\e H_1 \right) + \phi^\perp,\] we
need to estimate the two quantities
\begin{equation}
\label{lambdapartial}
 \int_{S_\e} \left(-\frac 2 \e |\n_{g_\e}
    u|^2 - 6 u_{k,\e} v_{k,\e} u^2\right) dV_{g_\e}; \qquad \quad
    \int_{S_\e} b u^2 dV_{g_\e}.
\end{equation}
Here the function $v_{k,\e}$ is defined as $v_{k,\e}=\frac{\pa
\tilde{u}_{k,\e}}{\pa \e}(\e \cdot)$, where $\tilde{u}_{k,\e}\,:\,
\e S_\e \to \R$ is the scaling $u_{k,\e}(x)=\tilde{u}_{k,\e}(\e x)$.
We claim that, normalizing $u$ with $\|u\|_{H^1(S_\e)} = 1$ (this
condition was required in Lemma \ref{l:charvarphi}), the following
estimates hold
\begin{equation}\label{lambdapartialexpansion}
\int_{S_\e} \left(-\frac 2 \e |\n_{g_\e}
    u|^2 - 6 u_{k,\e} v_{k,\e} u^2\right) dV_{g_\e} =
    \frac{c}{\e^{n-1}} \int_K b(\overline{y}) \var^2 + o(1);
\end{equation}
\begin{equation}\label{lambdapartialexpansion2}
 \int_{S_\e} b u^2 dV_{g_\e} =
    \frac{C_1}{\e^{n-1}} \int_K b(\overline{y}) \var^2 + o(1)
\end{equation}
as $\e \to 0$, where $c < 0$ and where $C_1$ is defined after
\eqref{eq:vvv}. This together with (\ref{eq:estdle}) would conclude
the proof of the lemma.

\

\noindent {\bf Proof of \eqref{lambdapartialexpansion} and
\eqref{lambdapartialexpansion2}.} First of all, recall that by our
normalization and by Lemma \ref{l:eigenlinop} we have that
$\|\phi^\perp\|_{H^1(S_\e)}=o(\e)$. Therefore, using the expansions
for $g_\e$ in Subsection \ref{ss:gb}, some integration by parts,
\eqref{eq:decPsi} and the estimates in Subsection \ref{ss:appsol}
one finds
\begin{eqnarray}\label{eq:estkato1}
  \nonumber&&\int_{S_\e} \left(-\frac 2 \e |\n_{g_\e}
    u|^2 - 6 u_{k,\e} v_{k,\e} u^2\right) dV_{g_\e} = \\
    && \frac2\e\int_{S_\e}\var (H'+\e H_1)
    \left[\var(H'''+\e H_1'')+\e\k \var(H''+\e H_1')+\e^2\D_K\var(H'+\e H_1) \right](1+\e \kappa \z)\\
  & +& 6\int_{S_\e}\var^2 (H+\e h_1)
    \left[\frac{\z-\Phi_0}{\e}H'+(\z-\Phi_0)h_1'\right](H'^2+2\e H'H_1)(1+\e \kappa \z)+o(1).\nonumber
\end{eqnarray}
Since the arguments of $H, H', h_1$ and $H_1$ are all translated by
$\Phi_0$ in $\z$, with the change of variables $s=\z-\Phi_0$ and
some elementary estimates (which use the exponential decay of $H',
h_1$ and $H_1$) we find
\begin{eqnarray}\label{eq:estkato12}
  \nonumber&&\int_{S_\e} \left(-\frac 2 \e |\n_{g_\e}
    u|^2 - 6 u_{k,\e} v_{k,\e} u^2\right) dV_{g_\e} = \frac 1\e
    \int_{K_\e \times \R} \var^2\left( 2H'H'''+6 s HH'^3
    \right) \\ & + & 4\int_{K_\e \times \R}\var^2 H_1
    \left(H'''+3 s HH'^2   \right) + 6 \int_{K_\e \times \R}\var^2 \bigg( s h_1 (H')^3+ s h'_1 H(H')^2
   \bigg) \\ & + & 2\int_{K_\e \times \R}\e \var \D_K\var H'^2+2\int_{K_\e \times \R}\k (s+\Phi_0) H'H'''\var^2
    +6 \int_{K_\e \times \R} s (s+\Phi_0) \k HH'^3\var^2+o(1). \nonumber
\end{eqnarray}
In the latter formula all the arguments now are simply in $s$, with
no more translation. Using equation \eqref{ode}, the oddness of $H$
together with some integration by parts, it is easy to see that the
term of order $\frac1\e$ in the above expression is identically
equal to zero. Let us consider now the terms of order $0$ which
involve $H_1$. Using the self-adjointness of the operator
$\mathcal{L}_0$ and the following elementary identity
$$
 \mathcal{L}_0\left( -\frac{1}{\sqrt{2}} s HH'  \right) = H'''+3s
 HH'^2,
$$
we can write that
\begin{equation}\label{eq:estkato2}
  \nonumber
  4\int_{K_\e \times \R}\var^2 H_1\left(H'''+3 s HH'^2   \right) =
   - \frac{4}{\sqrt{2}}\int_{K_\e \times \R}\var^2 s HH' \,\mathcal{L}_0 H_1.
\end{equation}
Similarly, the terms of order $0$ which involve $h_1$ can be written
(up to some integration by parts in the variable $s$) as
\begin{eqnarray}\label{eq:estkato3}
\nonumber
 6 \int_{K_\e \times \R}\var^2 \bigg( s h_1 (H')^3+ s h'_1 H(H')^2 \bigg)&=&
 - 6 \int_{K_\e \times \R}\var^2 \bigg( 2sH H'H''+ H(H')^2 \bigg)h_1\\
 & = & - 6 \int_{K_\e \times \R}\var^2 H(H')^2h_1+12\sqrt{2}\int_{K_\e \times \R}\var^2 s (H H')^2h_1\\
 & = & \frac{2}{\sqrt{2}} \int_{K_\e \times \R}\var^2   HH'\mathcal{L}_0(h_1)+\frac{24}{\sqrt{2}}
 \int_{K_\e \times \R}\var^2 s (H H')^2 h_1. \nonumber
\end{eqnarray}
Here again we have used the self-adjointness of the operator
$\mathcal{L}_0$ and the  identities
\begin{eqnarray*}
&&H''=-\sqrt{2} H H' \qquad \hbox{ and }\quad
   \mathcal{L}_0\left( \frac{2}{\sqrt{2}}  H H'  \right) = -6 H H'^2.
\end{eqnarray*}
Regrouping the above terms, and using the oddness of $H$ one finds
\begin{eqnarray}\label{eq:estkato4}
  \nonumber & & \int_{S_\e} \left(-\frac 2 \e |\n_{g_\e}
    u|^2 - 6 u_{k,\e} v_{k,\e} u^2\right) dV_{g_\e} \\ & = &
    \nonumber
    - \frac{4}{\sqrt{2}}\int_{K_\e \times \R}\var^2 s H H'
    \bigg(\mathcal{L}_0 H_1-6HH'h_1\bigg) \nonumber\\
    & + &\frac{2}{\sqrt{2}} \int_{K_\e \times \R}\var^2   HH'\mathcal{L}_0(h_1)
    + 2\int_{K_\e \times \R}\e \var \D_K\var H'^2\nonumber\\
    & + & 2\int_{K_\e \times \R}\k (s+\Phi_0) H'H'''\var^2
    + 6 \int_{K_\e \times \R}s (s+\Phi_0) \k HH'^3\var^2+o(1)\\
    & = & - \frac{4}{\sqrt{2}}\int_{K_\e \times \R}\var^2 s H H'
    \bigg(\mathcal{L}_0 H_1-6HH'h_1\bigg)
    +\frac{2}{\sqrt{2}} \int_{K_\e \times \R}\var^2   HH'\mathcal{L}_0(h_1) \nonumber \\
    &+&2\int_{K_\e \times \R}\e \var \D_K\var H'^2
    + 2 \int_{K_\e \times \R} \k \Phi_0 \var^2 (3 s HH'^3 + H'H''') +o(1).\nonumber
\end{eqnarray}
Now, using equations \eqref{eq:h1} and \eqref{eq:H1} we arrive at
\begin{eqnarray*}
  \nonumber & & \int_{S_\e} \left(-\frac 2 \e |\n_{g_\e}
    u|^2 - 6 u_{k,\e} v_{k,\e} u^2\right) dV_{g_\e}\
    \\ & = & \frac{4}{\sqrt{2}}\int_{K_\e \times \R}\var^2 s H H'
    \bigg(2 b(\e y) (s+\Phi_0) H H' + H'' \k(\e y) - \sqrt{2} b(\e y) H'\bigg)
   \\  & + &\frac{2}{\sqrt{2}} \int_{K_\e \times \R}\var^2   HH'\bigg( -\kappa(\e y) H'+
   (s+\Phi_0) b(\e y) (1-H^2) \bigg)\\  &+&2\int_{K_\e \times \R} \e \var \D_K\var H'^2
    + 2 \int_{K_\e \times \R} \k \Phi_0 \var^2 (3 s HH'^3 + H'H''')+o(1).\nonumber
\end{eqnarray*}
Using again the fact that $H$ is odd, by the vanishing of the last
integral (as one can easily check) we obtain
 \begin{eqnarray*}
  \nonumber\int_{S_\e} \left(-\frac 2 \e |\n_{g_\e}
    u|^2 - 6 u_{k,\e} v_{k,\e} u^2\right) dV_{g_\e}\
    & = & \frac{4}{\sqrt{2}}\int_{K_\e \times \R}\var^2 s H H'
    \bigg( 2 b(\e y) s H H'  -  \sqrt{2} b(\e y) H'\bigg)
   \\  & + &\frac{2}{\sqrt{2}} \int_{K_\e \times \R}\var^2   H H'
   s b(\e y) (1-H^2) +2\int_{K_\e \times \R} \e \var \D_K\var H'^2
   \nonumber \\ & + & o(1).\nonumber
\end{eqnarray*}
By an explicit computation of the integral we find
\begin{equation}\label{eq:ffin}
  \nonumber\int_{S_\e} \left(-\frac 2 \e |\n_{g_\e}
    u|^2 - 6 u_{k,\e} v_{k,\e} u^2\right) dV_{g_\e} = \left(
    \frac{8}{45} \pi^2 - \frac 23 \right) \frac{1}{\e^{n-1}} \int_{K} b \var^2 +
    \frac{4}{3} \frac{\e \sqrt 2}{\e^{n-1}} \int_{K} \var \D_K \var +
    o(1),
\end{equation}
so by Lemma \ref{l:charvarphi} and some easy estimates, we find
$$
  \nonumber\int_{S_\e} \left(-\frac 2 \e |\n_{g_\e}
    u|^2 - 6 u_{k,\e} v_{k,\e} u^2\right) dV_{g_\e} = \left(
    \frac{8}{45} \pi^2 - \frac{10}{3} \right)  \frac{1}{\e^{n-1}} \int_{K} b \var^2+o(1).
$$
Then we obtain \eqref{lambdapartialexpansion} taking
$c=\left(\frac{8}{45} \pi^2 - \frac{10}{3} \right)<0$. To prove
\eqref{lambdapartialexpansion2} it is sufficient to use Lemma
\ref{l:eigenlinop}, the estimates on $g_\e$ in Subsection
\ref{ss:gb} and the decay of $H_1$, see \eqref{eq:Psi}.
\end{pf}

\section{Proof of Theorem \ref{t:mmm}}\label{s:pf}

In this section we first prove the invertibility of the linearized
operator $L_\e$ using Lemma \ref{l:dereigenv} and choosing carefully
the parameter $\e$. Below, $H^2_0(S_\e)$ stands for the functions in
$H^2(S_\e)$ with null trace on $\pa S_\e$.

\begin{pro}\label{p:gap}
Let $k \geq 1$, let $u_{k,\e}$ be the approximate solution defined
in Proposition \ref{p:as}, and let $L_\e$ be the linearized operator
at $u_{k,\e}$, see \eqref{eq:linop}. Then there exist a sequence
$\e_j \to 0$ such that $L_{\e_j} : H^2_0(S_{\e_j}) \to
L^2(S_{\e_j})$ is invertible and its inverse $L^{-1}_{\e_j} :
L^2(S_{\e_j}) \to H_0^2(S_{\e_j})$ satisfies
$$
\left\| L_{\e_j}^{-1} \right\|_{\mathfrak{L}(L^2(S_\e);H^2_0(S_\e))} \leq C\e_j^{-\frac{n+1}{2}}, \qquad
\qquad \hbox{ for all } j \in \N.
$$
\end{pro}

\begin{pf}
The proof is similar in spirit to the one of Proposition 4.5 in \cite{mm2}.
 As we will see, in order to study the spectral gap of $L_\e$ it
suffices to find an asymptotic estimate on the number $N_\e$ of
positive eigenvalues of $L_\e$ and to apply then Lemma
\ref{l:dereigenv}. We denote by
$\tilde{\l}_{1,\e}\geq\tilde{\l}_{2,\e}\geq\cdots\geq
\tilde{\l}_{j,\e}\geq\cdots$ the eigenvalues of $L_\e$,  counted
with multiplicity. The $j$-th eigenvalue $\tilde{\l}_{j,\e}$ can be
estimated using the classical Courant-Fisher formula
\begin{equation}\label{eq:cf}
    \tilde{\l}_{j,\e} = \sup_{M \in M_j} \inf_{u \in M, u \neq 0}
    \frac{\int_{S_\e} u L_\e u dV_{g_\e}}{\int_{S_\e} b u^2 dV_{g_\e}};
    \qquad \quad \tilde{\l}_{j,\e} = \inf_{M \in M_{j-1}}
    \sup_{u \perp M, u \neq 0} \frac{\int_{S_\e} u L_\e u
    dV_{g_\e}}{\int_{S_\e} b u^2 dV_{g_\e}}.
\end{equation}
Here $M_l$ represents the family of $l$-dimensional subspaces of
$H^2_0(S_\e)$, and the symbol $\perp$ denotes orthogonality with
respect to the $L^2$ scalar product with weight $b$. Notice that the
$\inf$ and $\sup$ are reversed compared to \eqref{eq:ljcourant}
since the principal part of the operator has the opposite sign.

We can find a lower bound of $N_\e$ using the first formula in
\eqref{eq:cf}. Indeed, given a fixed $\d > 0$, let $j_\e$ be the
largest integer $j$ for which $\l_{j,\e} \geq \d \e$. From
\eqref{eq:weyllj} and \eqref{eq:lje} we find that
\begin{equation}\label{eq:asyje}
    j_\e \simeq \left( \frac{\sqrt{2} - \d}{C_{K,b} \e}
    \right)^{\frac{n-1}{2}} \qquad \quad \hbox{ as } \e \to 0.
\end{equation}
We can take a test function $\phi$ like $\var \Psi$ with
$\var=\sum_{l=1}^{j_\e} \a_l \var_l$. Actually, since we want to
work in the space $H^2_0(S_\e)$ we need to add a suitable cutoff
function in $\z$. However, by the exponential decay of $\Psi$, see
\eqref{eq:decPsi}, these will generate error terms exponentially
small in  $\e$. Therefore, for convenience of the exposition, we
will omit these corrections.

By \eqref{eq:LembbLe} we have that
$$
  L_\e \phi = (\e \sqrt2 b \var+\e^2 \D_{K}\var ) \Psi +\e^2 \var
  \D_{K} \Psi + 2 \e^2\langle \n_{K} \var, \n_{K}\Psi \rangle+\e^2 \tilde{R}_\e \var
+\e^{3}\hat{\mathfrak{L}}_{3,\e}(\var \Psi),
$$
where $\hat{\mathfrak{L}}_{3,\e}$ and $\tilde{R}_\e$ satisfy
respectively the estimates \eqref{eq:Lja} and \eqref{eq:Re}.
Reasoning as for \eqref{eq:esterror}, \eqref{eq:fethi} we find
\begin{equation}\label{eq:philephi}
\int_{S_\e}\phi\,L_\e\phi dV_{g_\e}\ge \frac{1}{\e^{n-1}}\sum_l
\left[ (1+o(1)) \l_{l,\e} + O(\e^2) \right] \a_l^2; \qquad \quad
\int_{S_\e} b \phi^2 dV_{g_\e} = \frac{1+o(1)}{\e^{n-1}} \sum_l
\a_l^2.
\end{equation}
Defining $M = span \{\var_l\Psi,\;\;l\le j_\e\}$, by the first
formula in \eqref{eq:cf} and our choice of $j_\e$ we have that
$$ \tilde{\l}_{j_\e,\e} \ge \inf_{u \in M, u \neq 0}
    \frac{\int_{S_\e} u L_\e u dV_{g_\e}}{\int_{S_\e} b u^2
    dV_{g_\e}} \geq 0, \qquad \quad \hbox{ as } \e \to 0.
$$
From \eqref{eq:asyje} and the last formula we then find the
following lower bound
\begin{equation}\label{eq:lbne}
  N_\e \geq (1 + o(1)) \left( \frac{\sqrt{2} - \d}{C_{K,b} \e}
    \right)^{\frac{n-1}{2}} \qquad \quad \hbox{ as } \e \to 0.
\end{equation}
A similar upper bound can be obtained using the second formula in
\eqref{eq:cf}: again, given a fixed $\d > 0$, let $\tilde{j}_\e$ be
the smallest integer $j$ for which $\l_{j,\e} \leq - \d \e$. Still
from \eqref{eq:weyllj} and \eqref{eq:lje} it follows that
\begin{equation}\label{eq:asyjte}
    \tilde{j}_\e \simeq \left( \frac{\sqrt{2} + \d}{C_{K,b} \e}
    \right)^{\frac{n-1}{2}} \qquad \quad \hbox{ as } \e \to 0.
\end{equation}
Now let $\phi\in H^2_0(S_\e)$ be an arbitrary function orthogonal to
$\tilde{M} := span \{\var_l\Psi,\;\;l\le \tilde{j}_\e-1\}$, and let
us write it in the form $\phi = \var \Psi + \phi^\perp$ with
$\phi^\perp$ as in Lemma \ref{l:eigenlinop}.

We write as before $\var = \sum_l \a_l\var_l$, and split it as sum
of the following two functions
$$
\ov{\var}_1= \sum_{l \leq \tilde{j}_\e-1} \a_l \var_l, \qquad \qquad
\ov{\var}_2= \sum_{l \geq \tilde{j}_\e} \a_l \var_l.
$$
Using the second formula in \eqref{eq:cf} we have that
\begin{equation}\label{eq:phijLephij}
  \tilde{\l}_{j,\e}\le \sup_{u \perp \tilde{M}, u \neq 0} \frac{\int_{S_\e} u L_\e u
    dV_{g_\e}}{\int_{S_\e} b u^2 dV_{g_\e}}.
\end{equation}
By the definition of $\phi^\perp$, the expansions of the metric
$g_\e$ in Subsection \ref{ss:gb} and \eqref{eq:decPsi} one finds
\begin{eqnarray}\label{eq:split}
  \nonumber \|\phi\|_{L^2(S_\e)}^2 & = & \|\var \Psi\|_{L^2(S_\e)}^2 +
  \|\phi^\perp\|_{L^2(S_\e)}^2 = (1 + o(1)) (\|\ov{\var}_1
  \Psi\|_{L^2(S_\e)}^2 +
   \|\ov{\var}_2\Psi\|_{L^2(S_\e)}^2)  + \|\phi^\perp\|_{L^2(S_\e)}^2 \\ & = & (1 +
  o(1)) \frac{1}{\e^{n-1}} \sum_{l \leq \tilde{j}_\e-1} \a_l^2 + (1 +
  o(1)) \frac{1}{\e^{n-1}} \sum_{l \geq \tilde{j}_\e} \a_l^2 +
  \|\phi^\perp\|_{L^2(S_\e)}^2 \qquad \hbox{ as } \e \to 0.
\end{eqnarray}
Using also the estimates in the proof of Lemma \ref{l:charvarphi},
multiplying $\phi$ by $b(\e y) \ov{\var}_1 \Psi$, using the
orthogonality to $\tilde{M}$ and integrating one can easily prove
that
\begin{equation}\label{eq:ovvphi1}
    \|\ov{\var}_1\Psi\|_{L^2(S_\e)} \leq
    o(1)\,\|\phi\|_{L^2(S_\e)}, \qquad \quad \hbox{ as } \e \to 0.
\end{equation}
This, together with \eqref{eq:LembbLe}, the fact that $|\l_{j,\e}|
\leq C \e$ for $j \leq \tilde{j}_\e - 1$ and some easy computations
imply
\begin{equation}\label{eq:ovv1leovv1}
  \int_{S_\e} (\ov{\var}_1\Psi) L_\e (\ov{\var}_1\Psi)
  dV_{g_\e} = o(\e) \|\phi\|_{L^2(S_\e)}^2.
\end{equation}
Now, by the self-adjointness of $L_\e$ we can write
\begin{eqnarray}\label{eq:tamara} \nonumber
  \int_{S_\e} \phi L_\e \phi dV_{g_\e} & = & \int_{S_\e} (\ov{\var}_1\Psi)
  L_\e (\ov{\var}_1\Psi) dV_{g_\e} + \int_{S_\e} (\ov{\var}_2\Psi) L_\e
  (\ov{\var}_2\Psi) dV_{g_\e} + \int_{S_\e} \phi^\perp
  L_\e \phi^\perp dV_{g_\e} \\
   & + & 2 \int_{S_\e} (\ov{\var}_1\Psi)
  L_\e (\ov{\var}_2\Psi) dV_{g_\e} + 2 \int_{S_\e} (\ov{\var}_1\Psi)
  L_\e \phi^\perp dV_{g_\e} + 2 \int_{S_\e} (\ov{\var}_2\Psi)
  L_\e \phi^\perp dV_{g_\e}.
\end{eqnarray}
Let us first estimate the last two terms: from \eqref{eq:phiperpLe}
we obtain
\begin{eqnarray}\label{eq:ovphi2phiperp}
 \int_{S_\e} (\ov{\var}_2\Psi)
  L_\e \phi^\perp dV_{g_\e} & \le & \frac{C\e^2}{\e^{\frac{n-1}{2}}}
  \|\phi^\perp\|_{H^1(S_\e)}\|\ov{\var}_2\|_{H^1(K)}\le
  \frac{C\e^2}{\e^{\frac{n-1}{2}}}\|\phi^\perp\|_{H^1(S_\e)}
  \left( \sum_{l\ge \tilde{j}_\e}(1+\l_{l})\a_l^2 \right)^{\frac 12}.
\end{eqnarray}
Similarly, using the fact that $|\l_l| \leq \frac C \e$ for $l \leq
\tilde{j}_\e - 1$ we have
$$
\|\ov{\var}_1\|_{H^1(K)}\le \left(\sum_{l \le \tilde{j}_\e - 1} \l_l
\a_l^2\right)^\frac12\le
\frac{1}{\e^\frac12}\|\ov{\var}_1\|_{L^2(K)},
$$
and hence
\begin{eqnarray}\label{eq:ovphi1phiperp}
 \int_{S_\e} (\ov{\var}_1\Psi)
  L_\e \phi^\perp dV_{g_\e} & \le & \frac{C\e^2}{\e^{\frac{n-1}{2}}}
  \|\phi^\perp\|_{H^1(S_\e)}\|\ov{\var}_1\|_{H^1(K)}\le
  \frac{C\e^\frac32}{\e^{\frac{n-1}{2}}}
  \|\phi^\perp\|_{H^1(S_\e)}\|\ov{\var}_1\|_{L^2(K)}.
\end{eqnarray}
On the other hand, a similar argument as for the first formula in \eqref{eq:philephi} yields
 \begin{equation}\label{eq:intovv2leovv2}
  \int_{S_\e} (\ov{\var}_2\Psi) L_\e (\ov{\var}_2\Psi) dV_{g_\e} =
  \frac{1}{\e^{n-1}}\sum_{l\ge \tilde{j}_\e}
  \left[(1+o(1))\l_{l,\e}+O(\e^{2})  \right]\a_l^2.
 \end{equation}
Moreover by the negative definiteness  of $L_\e$ on $\phi^\perp$,
see \eqref{eq:posphip}, we have that
\begin{equation}\label{eq:posphip2}
\int_{S_\e} \phi^\perp
  L_\e \phi^\perp dV_{g_\e} \leq - C^{-1} \| \phi^\perp\|^2_{H^1(S_\e)}
\end{equation}
for some fixed constant $C$. It remains to estimate the term
$\int_{S_\e} (\ov{\var}_1\Psi) L_\e (\ov{\var}_2\Psi) dV_{g_\e}$.
Using again \eqref{eq:LembbLe}, \eqref{eq:decPsi}, the fact that
$\int_K b \var_i \var_j = 0$ for $i \leq \tilde{j}_\e - 1$, $j \geq
\tilde{j}_\e$, and some integration by parts we get
\begin{eqnarray}\label{eq:mix}
  & & \frac{\e^{n-1}}{C} \int_{S_\e} (\ov{\var}_1\Psi)
  L_\e (\ov{\var}_2\Psi) dV_{g_\e} \leq \e^3 \int_K |\n \ov{\var}_1|
  \, |\n \ov{\var}_2| + \e^2 \int_K \left( |\n \ov{\var}_1| \, |\ov{\var}_2|
  + |\ov{\var}_1| \, |\n \ov{\var}_2| + |\ov{\var}_1|\,
  |\ov{\var}_2|\right) \\ \nonumber & \leq &
  \e^3 \left( \sum_{l\le \tilde{j}_\e-1} \l_{l}\a_l^2  \right)^{\frac12}
  \left( \sum_{l\ge \tilde{j}_\e} \l_{l}\a_l^2  \right)^{\frac12}+
 \e^2 \left( \sum_{l\le \tilde{j}_\e-1}\a_l^2  \right)^{\frac12}
  \left( \sum_{l\ge \tilde{j}_\e} \l_{l}\a_l^2  \right)^{\frac12} \\
   & + & \e^2 \left( \sum_{l\le \tilde{j}_\e-1} \l_{l}\a_l^2  \right)^{\frac12}
  \left( \sum_{l\ge \tilde{j}_\e} \a_l^2  \right)^{\frac12} +
  \e^2 \left( \sum_{l\le \tilde{j}_\e-1} \a_l^2  \right)^{\frac12}
  \left( \sum_{l\ge \tilde{j}_\e} \a_l^2  \right)^{\frac12}.
  \nonumber
  \end{eqnarray}
We claim next that the terms in the right-hand side of
\eqref{eq:tamara} can be combined yielding
\begin{equation}\label{eq:claim22}
  \int_{S_\e} \phi L_\e \phi dV_{g_\e} \leq \frac 1 C \left(
  \frac{1}{\e^{n-1}}\sum_{l\ge \tilde{j}_\e}
  \l_{l,\e} \a_l^2 - \| \phi^\perp\|^2_{H^1(S_\e)} \right) +
  o(\e) \|\phi\|_{L^2(S_\e)}^2.
\end{equation}

To prove this, we show that the main terms in \eqref{eq:tamara} are
the ones given in \eqref{eq:intovv2leovv2}, \eqref{eq:posphip2},
while all the others, listed in the left-hand sides of formulas
\eqref{eq:ovvphi1}, \eqref{eq:ovphi2phiperp},
\eqref{eq:ovphi1phiperp}, \eqref{eq:mix} can be absorbed into the
formers by the elementary inequality $|a b| \leq a^2 + b^2$. For
example, for any small constant $\b$ (independent of $\e$) we can
write
\begin{equation*}
  \frac{C\e^2}{\e^{\frac{n-1}{2}}}\|\phi^\perp\|_{H^1(S_\e)}
  \left( \sum_{l\ge \tilde{j}_\e}(1+\l_{l})\a_l^2 \right)^{\frac
  12} \leq C \b \|\phi^\perp\|_{H^1(S_\e)}^2 + \frac{C}{\b}
  \frac{1}{\e^{n-1}} \sum_{l\ge \tilde{j}_\e} \e^4
  (1+\l_{l})\a_l^2.
  \end{equation*}
Taking $\b$ sufficiently small, and noticing that $\e^4 \l_l \geq
\e^2 \l_{l,\e} + O(\e^3)$, from \eqref{eq:split} we deduce our claim
for this term. For the others, one reasons similarly, taking also
\eqref{eq:ovvphi1} and the choice of $\tilde{j}_\e$ into account.
Now \eqref{eq:split}, \eqref{eq:ovvphi1}, \eqref{eq:claim22} and
again the choice of $\tilde{j}_\e$ imply
$$
  \int_{S_\e} \phi L_\e \phi dV_{g_\e} \leq - \frac{\d \e}{C
  \e^{n-1}} \sum_{l\ge \tilde{j}_\e} \a_l^2 - \frac 1C
  \|\phi^\perp\|_{H^1(S_\e)}^2 \leq 0.
$$
Therefore, by \eqref{eq:asyjte} we find the following upper bound on
$N_\e$
\begin{equation*}\label{eq:ubne}
 N_\e \leq (1 + o(1))  \left( \frac{\sqrt{2} + \d}{C_{K,b} \e}
    \right)^{\frac{n-1}{2}} \qquad \quad \hbox{ as } \e \to 0.
\end{equation*}
Since $\d$ is arbitrary, the last estimate and \eqref{eq:lbne} imply
\begin{equation}\label{eq:nr+}
  N_{\e} \sim C_{1,K} \e^{-\frac{n-1}{2}} \qquad \qquad \hbox{ as } \e \to
  0,
\end{equation}
where we have set $C_{1,K}=\left(\frac{\sqrt{2}
}{C_{K,b}}\right)^{\frac{n-1}{2}}$.

Next, for $l \in \N$, let $\e_l = 2^{-l}$. From \eqref{eq:nr+} we
get
\begin{equation}\label{eq:diffe}
  N_{\e_{l+1}} - N_{{\e_{l}}} \sim C_{1,K}
  \left( 2^{(l+1)\frac{n-1}{2}} - 2^{l\frac{n-1}{2}} \right)
  = C_{1,K} (2^{\frac{n-1}{2}} - 1) \e_l^{-\frac{n-1}{2}}.
\end{equation}
On the other hand it follows from Lemma \ref{l:dereigenv} that the
eigenvalues of $L_\e$ which are bounded (in absolute value) by
$O(\e^2)$ are decreasing in $\e$. Equivalently, by the last
equation, the number of eigenvalues which become positive, when
$\e$ decreases from $\e_l$ to $\e_{l+1}$, is of order
$\e_l^{-\frac{n-1}{2}}$. Now we define
$$
\mathcal{A}_l = \left\{ \e \in (\e_{l+1}, \e_l) \; : \; ker L_\e
\neq \emptyset \right\}; \qquad \qquad \mathcal{B}_l = (\e_{l+1},
\e_l) \setminus \mathcal{A}_l.
$$
By \eqref{eq:diffe} and the monotonicity (in $\e$) of the {\em
small} eigenvalues, we deduce that card$(\mathcal{A}_l) < C
\e_l^{-\frac{n-1}{2}}$, and hence there exists an interval $(a_l,
b_l)$ such that
\begin{equation}\label{eq:albl}
  (a_l, b_l) \subseteq \mathcal{B}_l; \qquad \qquad |b_l - a_l| \geq C^{-1}
  \frac{\hbox{meas}(\mathcal{B}_l)}{\hbox{card}(\mathcal{A}_l)} \geq C^{-1} \e_l^{\frac{n+1}{2}}.
\end{equation}
From Lemma \ref{l:dereigenv} we deduce that $L_{{\frac{a_l+b_l}{2}}}$ is invertible and
$$
\left\|L_{{\frac{a_l+b_l}{2}}}^{-1}\right\|_{\mathfrak{L}(L^2(S_\e);H^2_0(S_\e))}
\leq \frac{C}{\e_l^{\frac{n+1}{2}}}.
$$
 This concludes the proof taking $\e_j = \frac{a_j+b_j}{2}$.
\end{pf}

\

\no Proposition \ref{p:gap} gives us a {\em localized} version of
the invertibility result we need. To have a global one in the whole
domain $\O_\e$, we define a smooth cutoff function $\chi_\e$ by
$$
\left\{
  \begin{array}{ll}
    \chi_\e(t)= 1, & \hbox{for } t \leq \frac 12 \e^{-\g} \\[2mm]
    \chi_\e(t)= 0, & \hbox{for }t \geq \frac 34 \e^{-\g}\\[2mm]
    \left|\chi'_\e\right|\le C\e^{\g} & \hbox{and } \left|\chi''_\e\right|\le C\e^{2\g}.
  \end{array}
\right.
$$
We next set
\begin{equation}\label{eq:hatuke}
\hat{u}_{k,\e}(y,\z) = \left\{
          \begin{array}{ll}
          1 + \chi_\e(\z) \left(u_{k,\e}(y,\z)-1\right), & \hbox{in } \O_+ \\
          -1 -  \chi_\e(\z) \left(u_{k,\e}(y,\z)+1\right), & \hbox{in } \R^n\setminus
          \O_+.
         \end{array}
         \right.
\end{equation}
Now by \eqref{asyw} and \eqref{eq:decuke} we know that $|u_{k,\e}|$
is exponentially close to $1$ and its derivative are exponentially
small for $|\z|$ sufficiently large. This and \eqref{eq:deceqas}
imply that $\|\mathfrak{S}_\e(\hat{u}_{k,\e})\|_{L^2(\O_\e)} \leq C
\e^{k+1-\frac{n-1}{2}}$ and
$\|\mathfrak{S}_\e(\hat{u}_{k,\e})\|_{L^\infty(\O_\e)} \leq C
\e^{k+1}$, see \eqref{eq:mfse}, where $C$ depends only on $K$, $b$
and $k$.

Let us denote by $\hat{L}_\e$ the linearized operator at
$\hat{u}_{k,\e}$ in $\O_\e$. Given a smooth positive extension
$\hat{b}$ of $b$ to $\ov{\O}$, we consider next the eigenvalue
problem
\begin{equation}\label{eq:eigOe}
    \left\{
      \begin{array}{ll}
        \hat{L}_\e u = \l \hat{b}(\e \cdot) u & \hbox{ in } \O_\e; \\
        \frac{\pa u}{\pa \nu} = 0 & \hbox{ on } \pa \O_\e,
      \end{array}
    \right.
\end{equation}
and we denote its eigenvalues by $(\hat{\l}_{j,\e})_j$, counted in
decreasing order with their multiplicity.

By \eqref{eq:linop}, asymptotically away from $K_\e$, the
eigenfunctions $u$ satisfy
$$
  \D u - (2 \pm 2 a(\e x) - \l \hat{b}) u = 0 \qquad \hbox{ in }
(\O_{\pm})_\e,
$$
where $(\O_{\pm})_\e$ stands for the $\frac 1 \e$ dilation of
$\O_{\pm}$.  Since we are assuming $a \in (-1,1)$ in $\ov{\O}$, if
$\l$ is bounded from below by $-\e$, the coefficient of $u$ in the
above equation is negative. Hence, reasoning as in \cite{mm1}, Lemma
5.1, one can prove that $u$ has an exponential decay away from
$K_\e$.

Moreover, an argument based on the Courant-Fisher method, see
Proposition 5.6 in \cite{mm1}, shows that there exists a constant
$C$ depending only on $\O$, $K$, $a$ and $\g$ such that
$$
  |\hat{\l}_{j,\e} - \tilde{\l}_{j,\e}| \leq C e^{-\frac{C}{\e}}
  \qquad \hbox{ provided } \hat{\l}_{j,\e} \geq - \e \hbox{ or }
  \tilde{\l}_{j,\e} \geq - \e.
$$
Here $\tilde{\l}_{j,\e}$ are the eigenvalues of $L_\e$ in $S_\e$,
see the proof of Proposition \ref{p:gap}.

This and Proposition \ref{p:gap} allow us to prove the following
result which guarantee the invertibility of the linearized operator for the range
 of the parameter $\e$ constructed above.

\begin{corollary}\label{c:gap}
Fix $k \in \N$ and let $\hat{u}_{k,\e}$, $\hat{L}_\e$ be as above.
Define $H^2_\nu(\O_\e)$ to be the subset of $H^2(\O_\e)$ consisting
of the functions with zero normal derivative at $\pa \O_\e$. Then
for a suitable sequence $\e_j \to 0$, the operator $\hat{L}_{\e_j} :
H^2_\nu(\O_{\e_j}) \to L^2(\O_{\e_j})$ is invertible and the inverse
operator satisfies $\left\| \hat{L}_{\e_j}^{-1}
\right\|_{\mathfrak{L}(L^2(\O_{\e_j});H^2_\nu(\O_{\e_j}))} \leq
\frac{C}{\e_j^{\frac{n+1}{2}}}$, for all $j \in \N$.
\end{corollary}

\

\no Using the above results, we are in position to prove our main result, Theorem
\ref{t:mmm}.

\

\begin{pfn} {\sc of Theorem \ref{t:mmm}} Let $(\e_j)_j$ be as in Corollary \ref{c:gap}.
We look for a solution $u_\e$ of the equation
$\mathfrak{S}_\e(u)=0$ of the form
$$
u_\e = \hat{u}_{k,\e} + w, \qquad w \in H^2_\nu(\O_\e).
$$
For $\e=\e_j$, define the function $\hat{F}_\e: H^2_\nu(\O_\e) \cap
L^\infty(\O_\e)\to H^2_\nu(\O_\e) \cap L^\infty(\O_\e)$ by
\begin{equation}\label{eq:hatfe}
 \hat{F}_\e(w) := - \hat{L}_\e^{-1} \left[ \mathfrak{S}_\e(\hat{u}_{k,\e})
 - \left(3 \hat{u}_{k,\e}-a\right) w^2 - w^3 \right].
\end{equation}
We have that
\begin{equation}\label{eq:taut}
  \mathfrak{S}_\e(\hat{u}_{k,\e} + w) = 0 \qquad \Longleftrightarrow \qquad  \hat{F}_\e(w)=w.
\end{equation}
We want to prove that $\hat{F}_\e$ is a contraction in some closed
ball of $H^2_\nu(\O_\e) \cap L^\infty(\O_\e)$. We first define the
norm $|||\cdot|||$ as $|||w||| = \|w\|_{H^2_\nu(\O_\e)} +
\|w\|_{L^\infty(\O_\e)}$. Then, for $r > 0$, we introduce the set
$$
  \mathcal{B}_r = \left\{ w \in H^2_\nu(\O_\e) \cap L^\infty(\O_\e) \;
  : \; |||w||| \leq r \right\}.
$$
Applying a standard elliptic regularity theorem and using
Corollary \ref{c:gap} one can prove that there exists positive
constants $C$ (depending on $\O$, $K$ and $a$) and $d$ (depending on
the dimension $n$) such that
\begin{eqnarray}
 |||\hat{F}_\e(w)||| &\leq& C \e^{-d} \left( \e^{k+1-\frac{n-1}{2}} +
  |||w|||^2 \right);\\
  |||\hat{F}_\e(w_1) - \hat{F}_\e(w_2)||| &\leq& C \e^{-d} \left(
  |||w_1||| + |||w_2||| \right) (|||w_1 - w_2|||),\nonumber
\end{eqnarray}
for $\e = \e_j$ and $w, w_1, w_2 \in H^2_\nu(\O_\e) \cap
L^\infty(\O_\e)$. Now setting $r = \e^l$, we can choose first $k$
sufficiently large, depending on $d$, $n$ and then $l$ depending on
$d$ and $k$ so that $\hat{F}_\e$ is a contraction in the ball
$B_{r}$ for $\e=\e_j$ sufficiently small. A solution of
\eqref{eq:taut} can be then found using the contraction mapping
theorem and its properties follows from the construction of
$u_{k,\e}$. This concludes the proof.
\end{pfn}

\begin{center}
{\bf Acknowledgments}
\end{center}
F. M and A. M. are supported by MURST, within the PRIN 2006 {\em
Variational Methods and Nonlinear Differential Equations}. The
research of J. W. is partially  supported by an Earmarked Grant from
RGC   of Hong Kong.

\

\end{document}